\documentclass{amsart}

\usepackage{amssymb}
\usepackage{amsmath}
\usepackage{array}
\usepackage{graphicx}
\usepackage{multirow}
\usepackage{xcolor}

\usepackage{mathtools}

\newcommand{\Pf}{{\em Proof}. }
\newcommand{\EPf}{\hfill$\square$}

\newcommand{\Lg}{\mbox{$\mathfrak g$}}

\newcommand{\Lk}{\mbox{$\mathfrak k$}}
\newcommand{\Lp}{\mbox{$\mathfrak p$}}
\newcommand{\La}{\mbox{$\mathfrak a$}}

\usepackage{accents}

\newcommand{\vardbtilde}[1]{\tilde{\raisebox{0pt}[0.85\height]{$\tilde{#1}$}}}

\newcommand\liegr{\sf}

\newcommand{\SU}[1]{\mbox{${\liegr SU}(#1)$}}

\newcommand{\U}[1]{\mbox{${\liegr U}(#1)$}}
\newcommand{\SP}[1]{\mbox{${\liegr Sp}(#1)$}}
\newcommand{\SO}[1]{\mbox{${\liegr SO}(#1)$}}

\newcommand{\OG}[1]{\mbox{${\liegr O}(#1)$}}
\newcommand{\Spin}[1]{\mbox{${\liegr Spin}(#1)$}}
\newcommand{\G}{\mbox{${\liegr G}_2$}}
\newcommand{\F}{\mbox{${\liegr F}_4$}}
\newcommand{\E}[1]{\mbox{${\liegr E}_{#1}$}}

\newcommand\fieldsetc{\mathbb}

\newcommand{\Z}{\fieldsetc{Z}}
\newcommand{\R}{\fieldsetc{R}}
\newcommand{\C}{\fieldsetc{C}}
\newcommand{\Q}{\fieldsetc{H}}
\newcommand{\Ca}{\fieldsetc{O}}
\newcommand{\K}{\fieldsetc{K}}

\newtheorem{thm}{Theorem}[section]
\newtheorem{cor}[thm]{Corollary}
\newtheorem{prop}[thm]{Proposition}
\newtheorem{lem}[thm]{Lemma}

\theoremstyle{remark}
\newtheorem{rmk}[thm]{Remark}

\DeclareMathOperator{\rk}{rk}
\DeclareMathOperator{\nsf}{nsf}

\newcommand{\ft}{\footnote}

\title{Actions on positively curved manifolds
  and\\ boundary in the orbit space}
\author{Claudio Gorodski}
\thanks{The first author acknowledges partial financial
  support from  CNPq (grant 302882/2017-0) and FAPESP (grant 16/23746-6).}

\address{Instituto de Matem\'atica e Estat\'\i stica, Universidade de 
  S\~ao Paulo, Rua do Mat\~ao, 1010, S\~ao Paulo, SP 05508-090, Brazil}
\email{gorodski@ime.usp.br}

\author{Andreas Kollross}
\address{Universität Stuttgart, Institut für Geometrie und Topologie,
Pfaffenwaldring~57, 70569~Stuttgart, Germany}

\email{kollross@mathematik.uni-stuttgart.de}

\author{Burkhard Wilking}

\address{Mathematisches Institut, M\"unster University, Einsteinstrasse~62,
  48149~M\"unster, Germany}
\email{wilking@uni-muenster.de}

\date{\today}

\subjclass[2010]{57S15 53C21 53C35}

\begin{document}

\maketitle

\begin{abstract}
  We study isometric actions of compact Lie groups on
  complete orientable positively curved $n$-manifolds whose orbit spaces
  have non-empty boundary in the sense of Alexandrov geometry.
  In particular,
  we classify quotients of the unit sphere by actions of compact simple Lie groups
  with non-empty boundary.
 We deduce from this the list of representations
 of compact simple Lie groups that admit non-trivial reductions.
 As a tool of special interest, we introduce a new geometric invariant
 of a compact symmetric space, namely, the minimal number of points
 in a ``spanning set'' of the space.
\end{abstract}

\section{Introduction}

\subsection{General observations}

For an isometric action of a compact Lie group $G$ on a
complete Riemannian manifold $M$ with orbit space $X=M/G$ stratified
by orbit types, the boundary of $X$ consists of the most important singular
strata of $X$; here the \emph{boundary}~$\partial X$ is defined as the closure
of the union of all strata of codimension one of~$X$.
In case $M$ is positively curved, this notion of boundary coincides with the
boundary of $X$ as an Alexandrov space and has a bearing on the geometry and
topology of $X$. For instance it is easy to see
that $\partial X$ is non-empty if and only
if $X$ is contractible (for the 'only if' part one uses the fact that
the distance to the boundary is a strictly concave function hence admits a
unique point of maximum, a ``soul point''; the 'if' part follows from the fact
the Alexander-Spanier $\Z_2$-cohomology in top degree of $X$ is
non-trivial if $\partial X=\varnothing$~\cite[Lemma~1]{GP}).
In general, the boundary plays an important role in some proofs
in the literature; see e.g.~main results in~\cite{S}, or
\cite[Theorem~1.4]{AR} and~\cite[\S5.3]{GL}.

\subsection{The case of quotients of the sphere}
It follows from the slice theorem that the
presence of boundary is a local condition,
in the sense that $X=M/G$ has non-empty boundary if and only if there
exists a point $p\in M$ such that the slice representation of
the isotropy group $G_p$ on the normal space $\nu_p(Gp)$ to the orbit~$Gp$
has orbit space with non-empty boundary. The orbit space of
an orthogonal representation is a metric cone over the
orbit space of the corresponding unit sphere, so also the
boundary of the former is a metric cone over the boundary of the latter.
These remarks show that
the special case of quotients of the unit sphere with non-empty
boundary plays a distinguished role.

In fact, as a main consequence of our methods, we deduce a rather
simple criterion for the existence of boundary
for quotients of spheres (or more generally, positively curved manifolds)
by \emph{simple} groups.

\begin{thm}\label{simple}
  Let $G$ be a compact connected simple Lie group. Then there is an
  explicit, positive integer $\mathcal L_G$, depending only on the local
  isomorphism class of $G$, such that: For every effective and isometric
  action of $G$ on a connected complete orientable
  Riemannian manifold $M$ of positive sectional curvature,
  if $\dim M\geq\mathcal L_G$ and the orbit-space has non-empty boundary, then
  the $G$-fixed point set $M^G\neq\varnothing$ and
  $\dim M^G\geq\dim M-\mathcal L_G$.
\end{thm}

The number $\mathcal L_G$ is easy to determine (cf.~Table~4 for its values)
and has geometric meaning, namely
\begin{equation}\label{lg}
 \mathcal L_G:= \max_K\{\ell_{G/K}(4+\dim G/K)\}, 
\end{equation}
where $K$ runs through all symmetric subgroups of $G$ with
maximal rank, and $\ell_{G/K}$ is defined as the minimum number
$\ell$ such that there exist
$\ell$ points in $G/K$ not contained in any proper closed
connected totally geodesic submanifold (cf.~section~\ref{tg}).

The number $\ell_{G/K}$ is a natural,
geometric invariant of a compact symmetric space $G/K$,
which is related to the minimum number of
involutions of $G$ necessary to topologically generate the group
(Proposition~\ref{dense}), and to the
minimum number of generic points of $G/K$ which are not
simultaneously fixed by a
non-identity isometry in $G$ (Proposition~\ref{h}).
In a sense, it is the minimum number of
points ``spanning'' $G/K$, and loosely alludes to the concept of linearly
independence in Linear Algebra. For instance, for the sphere we have
$\ell_{S^n}=n+1$. However,
the case of rank one symmetric spaces and Grassmannians
turns out to be special, as $\ell_{G/K}=3$ for the other spaces
that we compute (cf.~Theorem~\ref{ell}).

Applying Theorem~\ref{simple} to orthogonal actions on unit
spheres yields that a representation of a compact connected
simple Lie group $G$ on an Euclidean space $V$ that has no trivial
components can have orbit space
with non-empty boundary only if $\dim V\leq\mathcal L_G$.
We obtain a classification of
such representations by combining this remark with
a result about reducible representations (Corollary~\ref{cor:reducible}).

\begin{thm}\label{classif}
  The representations $V$ of
compact connected simple Lie groups $G$ with non-empty boundary
in the orbit space are listed in Tables~1 and~2,
up to a trivial component and up to an outer automorphism.
In the irreducible case (Table~1),
we also indicate the kernel of the representation in those cases
in which it is non-trivial, the effective
principal isotropy group,
and whether the representation is polar, toric
or quaternion-toric (we recall these concepts in subsection~\ref{sec:cplx}).

%{\renewcommand{\arraystretch}{1.3}
\begin{table}[t]
\rm
\[ \begin{array}{|c|c|c|c|c|}
  \hline
  G & \textsl{Kernel} & V & \textsl{Property} & \textsl{Effective p.i.g.} \\
  \hline
  \SU2 &\mbox{---}& \C^2 & \textrm{polar} & 1 \\
  \hline
   \multirow{2}{*}{\SO3}
  &\multirow{2}{*}{\mbox{---}}& \R^3 &\multirow{2}{*}{polar}&\sf T^1\\
  &&\mathrm{S}^2_0\R^3=\R^5 && \sf \Z_2^2\\
  \hline
  \multirow{3}{*}{\shortstack{\SU n\\ ($n\geq3$)}} &\mbox{---}&\C^n & \multirow{2}{*}{\textrm{polar}}&\SU{n-1}\\
 &\Z_n&\mathrm{Ad}&&\sf T^{n-1} \\ \cline{4-4}
  &\textrm{$\{\pm1\}$ if $n$ is even} &\mathrm{S}^2\C^n& \textrm{toric}&\Z_2^{n-1}/\mathrm{ker}\\
  \hline
    \multirow{2}{*}{\shortstack{\SU n\\ ($n\geq5$)}}
   &\multirow{2}{*}{\textrm{$\{\pm1\}$ if $n$ is even}}
   &\multirow{2}{*}{$\Lambda^2\C^n$}
   &\parbox[t]{.8in}{polar ($n$~odd)}
   &\multirow{2}{*}{$\SU2^{\lfloor\frac n2\rfloor}/\mathrm{ker}$}\\
   &&&\parbox[t]{.8in}{toric ($n$~even)}&\\
  \hline
  \SU6 &\mbox{---}& \Lambda^3\C^6=\Q^{10} & \textrm{q-toric}&\sf T^2\\ \hline
  \SU8& \Z_4 & [\Lambda^4\C^8]_{\mathbb R} & \textrm{polar} & \Z_2^7 \\ \hline
  \multirow{3}{*}{\shortstack{\SO n\\ ($n\geq5$)}}  &\mbox{---}& \R^n &
  \multirow{3}{*}{polar}&\Spin{n-1}\\ \cline{2-2}
  &\multirow{2}{*}{\textrm{$\{\pm1\}$ if $n$ is even}}&\Lambda^2\R^n=\mathrm{Ad}&&\sf T^{\lfloor\frac n2\rfloor} \\
  &&\mathrm{S}^2_0\R^n&&\Z_2^{n-1}/\mathrm{ker}\\  \hline
  \Spin7&\mbox{---}&\R^8\;\textrm{(spin)}& \textrm{polar} &\G \\ \hline
  \Spin8 &\Z_2& \R^8_\pm\; \textrm{(half-spin)} &\textrm{polar} & \mathsf{Spin}'(7) \\  \hline
  \Spin9 & \mbox{---} &\R^{16}\;\textrm{(spin)} & \textrm{polar} & \Spin7 \\ \hline
  \Spin{10} & \mbox{---} &\C^{16}_\pm\; \textrm{(half-spin)} & \textrm{polar} & \SU4 \\ \hline
  \Spin{11} & \mbox{---}& \Q^{16}\;\textrm{(spin)} & \mbox{---} & 1 \\ \hline
  \Spin{12} & \Z_2 & \Q^{16}_\pm\;\textrm{(half-spin)} & \textrm{q-toric} & \SP1^3\\ \hline
    \Spin{16} & \Z_2 & \R^{128}_\pm\;\textrm{(half-spin)} & \textrm{polar} & \Z_2^8\\ \hline
  \multirow{3}{*}{\shortstack{\SP n\\ ($n\geq3$)}} &\mbox{---}& \C^{2n}=\Q^n & \multirow{3}{*}{polar}&\SP{n-1} \\  \cline{2-2}
  &\multirow{2}{*}{$\{\pm1\}$}&[\mathrm{S}^2\C^{2n}]_{\mathbb R}=\mathrm{Ad}&&\sf T^n \\
  &&[\Lambda^2_0\C^{2n}]_{\mathbb R} &&\SP1^n/\{\pm1\} \\
  \hline
  \SP3 &\mbox{---}& \Lambda^3_0\C^6=\Q^7 &\textrm{q-toric}&\Z_2^2 \\ \hline
  \SP4 & \{\pm1\} & [\Lambda^4_0\C^8]_{\mathbb R} & \textrm{polar}&\Z_2^6 \\ \hline
  \multirow{2}{*}{\G} &\multirow{2}{*}{\mbox{---}}& \R^7 & \multirow{2}{*}{polar} &\SU3 \\
  &&\mathrm{Ad} &&\sf T^2\\ \hline
  \multirow{2}{*}{\F} &\multirow{2}{*}{\mbox{---}}& \R^{26}
  &\multirow{2}{*}{polar}&\Spin8\\
  &&\mathrm{Ad} &&\sf T^4\\ \hline
  \E6 &\mbox{---}& \C^{27} & \textrm{toric} & \Spin8\\ \hline
  \E6& \Z_3 &\mathrm{Ad} & \textrm{polar} & \sf T^6\\ \hline
  \E7 &\mbox{---}& \Q^{28} & \textrm{q-toric} & \Spin8\\ \hline
    \E7 & \Z_2 &\mathrm{Ad} &  \textrm{polar} & \sf T^7 \\\hline
  \E8 &\mbox{---}& \mathrm{Ad} & \textrm{polar} & \sf T^8\\
   \hline
\end{array} \]
\smallskip
\begin{center}
  \sc Table~1: Irreducible representations of compact simple
  Lie groups with non-empty boundary in the orbit space.
\end{center}
\end{table}
%}

\begin{table}[t]
\[ \begin{array}{|c|c|c|}
  \hline
  \multirow{2}{*}{\SU n} & k\,\C^n & 2\leq k\leq n-1 \\
  & \C^n\oplus \Lambda^2\C^n  & n\geq5 \\
  \hline
  \multirow{2}{*}{\SU4} & k\,\R^6\oplus\ell\,\C^4 &  2\leq k+\ell\leq 3\\
  & \R^6\oplus\mathrm{Ad} & - \\
  \hline
  \multirow{2}{*}{\SO n} & k\,\R^n & 2\leq k\leq n-1\\
  & \R^n\oplus\mathrm{Ad} & n\geq5 \\
  \hline
  \SP2 & \Q^2\oplus\R^5 & -\\
  \hline
  \Spin7 & k\,\R^7\oplus\ell\,\R^8 & 2\leq k+\ell\leq 4 \\
  \hline
  \Spin8 & k\,\R^8\oplus\ell\,\R^8_+\oplus m\,\R^8_- & 2\leq k+\ell+m\leq5\\
  \hline
  \multirow{3}{*}{\Spin9} & k\, \R^{16} & 2\leq k\leq3\\
  & \R^{16}\oplus k\,\R^9 & 1\leq k\leq4\\
  &2\R^{16}\oplus k\,\R^9 & 0\leq k\leq2\\
  \hline
  \Spin{10} & \C^{16}\oplus k\,\R^{10} & 1\leq k\leq 3\\
  \hline
  \Spin{12} & \Q^{16}\oplus\R^{12} & -\\
  \hline
  \multirow{2}{*}{\SP n} & k\,\C^{2n} & 2\leq k\leq n\\
  & \C^{2n} \oplus[\Lambda^2_0\C^{2n}]_{\mathbb R} & n\geq3 \\ \hline
  \SP3 & 2\,[\Lambda^2_0\C^6]_{\mathbb R} & - \\ \hline
  \G & k\,\R^7 & 2\leq k\leq3\\
  \hline
  \F & 2\,\R^{26} & -\\
  \hline
\end{array} \]
\smallskip
\begin{center}
  \sc Table~2: Reducible representations of compact simple
  Lie groups with non-empty boundary in the orbit space.
\end{center}
\footnotesize{\rm In case of $\Spin8$,
  the prime in $\mathsf{Spin}'(7)$ refers to a nonstandard $\Spin7$-subgroup;
  in case of $\Spin n$, $\mathrm{S}^2_0\R^n=\mathrm{S}^2\R^n\ominus\R$;
  in case of $\SP n$,
  $\Lambda^k_0\C^{2n}=\Lambda^k\C^{2n}\ominus\Lambda^{k-2}\C^{2n}$;
  and $[V]_{\mathbb R}$ denotes a real form of~$V$.}
\end{table}
\end{thm}

To exemplify the usefulness of the remark about the
existence of boundary being a local property, we give the
following result. The special thing about the groups listed
in the statement of the next proposition is that
according to Theorem~\ref{classif}
they are simple Lie groups for which a given representation
has non-empty
boundary in the orbit space if and only if it is polar.

\begin{cor}\label{some-simple}
  Let $G$ be one of the following simple Lie groups:
  \begin{gather*}
    \SU2,\ \SU n/\Z_n\ \textrm{($n\geq3$)},\ \SU8/\Z_4,\
  \SO n/\{\pm1\}\ \textrm{($n\geq6$ even)},\\ \mathsf{SO}'(16),\
  \SP n/\{\pm1\}\ \textrm{($n\geq4$)},\ \E6/\Z_3,\ \E7/\Z_2,\ \E8
  \end{gather*}
  ($\mathsf{SO}'(16)$ denotes a group isomorphic to the
  image of $\Spin{16}$ under a
  half-spin representation). Consider an \textbf{effective} isometric action of $G$
on a connected simply-connected
compact Riemannian manifold $M$ of positive
sectional curvature and dimension $n>\mathcal L_G$ (see Table~4 for the
explicit values of $\mathcal L_G$).
Then the orbit space $X=M/G$ has non-empty boundary
if and only if the action is polar; further, in this case
$M$ is equivariantly diffeomorphic a compact rank one symmetric space
with a linearly induced action.
\end{cor}

\subsection{The complexity of orbit spaces}\label{sec:cplx}

Our results also have a bearing on understanding the
``complexity'' of quotients of the unit sphere.
In the case of
orthogonal representations of a compact Lie group on
vector spaces (or more generally, isometric actions on
positively curved manifolds), the following criteria have been used to
describe representations whose geometry is not too complicated, namely:
\begin{enumerate}
  \item[(i)] The principal isotropy group is non-trivial.
  \item[(ii)] There exists a non-trivial reduction, that is, a
    representation of a group with smaller dimension and isometric orbit space.
  \item[(iii)] The cohomogeneity, or codimension of the principal
    orbits, is ``low''.
%  \item[(iv)] There are ``few'' orbit types~\cite{HS}.
\end{enumerate}
It is known that (i) implies (ii)~\cite{St}, and
(ii) implies having non-empty boundary~\cite[Proposition~5.2]{GL}.
Indeed in case (i), the number of faces of the
boundary of the orbit space of an isometric action on a
positively curved manifold controls the number of simple factors
and the dimension of the center of the principal isotropy
group~\cite[Corollary~12.1]{Wi2}; here a \emph{face} is defined as
the closure of a component of a codimension one stratum.
We see a posteriori that to some extent (iii)
is also related to having non-empty boundary~\cite{HL}.
Representations with non-trivial principal isotropy group have been
partially classified in~\cite{HH} (however, note that the spin representation
of $\Spin{14}$ listed in Table~A therein
indeed has trivial principal
isotropy group; cf.~\cite[Remark~3.2]{Go}),
and the systematic study of representations with
non-trivial reductions (beyond polar representations)
has been initiated in~\cite{GL}.

Recall
that a representation is called \emph{polar}
if it admits a reduction to a representation of a finite group,
and it is called \emph{toric} (resp.~\emph{quaternion-toric})
if it is non-polar and it admits a reduction to a
representation of a group whose identity component
is Abelian (resp.~is isomorphic to~$\SP1^k$ for some $k>0$).
These classes are mostly related to the isotropy representations
of symmetric spaces. Polar representations are
classified in~\cite{D} (see also~\cite{B}).
Toric irreducible representations are classified in~\cite{GL2}
(see also~\cite{Pa} for some partial results in the reducible case).
Quaternion-toric irreducible representations are classified
in~\cite{GG}.

As another corollary to Theorem~\ref{classif}, we deduce:

\begin{cor}\label{cor:red}
  An irreducible representation of a compact connected simple Lie group
  admits a non-trivial reduction if and only if it is
  polar, toric or q-toric.
\end{cor}

Up to orbit-equivalence,
the representations 
in Corollary~\ref{cor:red} also coincide
with the representations of compact connected simple Lie groups
with non-trivial principal isotropy group~\cite[ch.~I, \S~2]{HH}.
Further, their minimal reductions are obtained from the fixed point
set of a principal isotropy group, after possibly enlarging
the group to an orbit-equivalent action. The (complexification
  of the) isometry between
  the orbit spaces given by this kind of reduction was shown to be
an isomorphism of affine algebraic varieties in~\cite{luna-richardson};
in particular, it is a diffeomorphism in the sense of~\cite{S}. In this
sense, Corollary~\ref{cor:red} can also be seen as a small step
toward proving the conjecture that 
a version of the Myers-Steenrod theorem holds for orbit spaces,
namely, that the smooth structure
is determined by the metric structure
(see~\cite[\S1.1, 1.2, 1.3]{alexandrino-lytchak}
and~\cite[\S1]{AR}).

\subsection{Quaternionic representations}

The following result came out of discussions of the first
named author with Ricardo Mendes.
It implies that the identity component of the isometry group
of the orbit space of an irreducible representation of
quaternionic type with cohomogeneity at least two is isomorphic
to $\SP1$ or $\SO3$ (compare~\cite{mendes}). 

\begin{cor}\label{quat}
  Let $\rho:G\to\OG V$ be
  an irreducible representation of quaternionic
type of a compact connected Lie group $G$ with cohomogeneity 
$c(\rho)\geq2$.
Consider the natural enlargement
$\hat\rho:\hat G\to\OG V$, where $\hat G=G\times\SP1$.
Then the cohomogeneities of these representations satisfy
\[ c(\rho)=c(\hat\rho)+3. \]
In particular, $\hat\rho$ is not
orbit-equivalent to~$\rho$.
\end{cor}

\subsection{Dimension estimate}

After a presentation of our applications, we have now come
to the rather technical statement of our most general
main result, although in the present paper
we have not had the opportunity of applying it in its full force.
It is a general estimate on the
dimension of a positively curved
manifold on which a Lie group acts with orbit space with
non-empty boundary. The normal subgroup $N$
in Theorem~\ref{main} contains all the information about the
boundary of $X$ and has a fixed point;
its existence is an act of balance between condition (a)
that restricts the largeness of $N$, and
condition~(c) that restricts its smallness.
Note that in case $G$ is simple, the theorem is just
saying that $G$ has a fixed point.

\begin{thm}\label{main}
  Let $G$ be a compact connected Lie group acting
  isometrically and effectively on a connected complete orientable
$n$-manifold $M$ of positive sectional curvature.
  Assume that $X=M/G$ has non-empty boundary and
  \begin{equation}\label{ineq}
    n>\alpha_G + \beta_G
  \end{equation}
  where
  \[ \alpha_G=2\dim G_{ss}+8\rk G_{ss} + 4\nsf G_{ss}\quad\mbox{and}\quad
  \beta_G=2\dim Z(G); \]
  here $Z(G)$ denotes the center of $G$, $G_{ss}=G/Z(G)$ its
  semisimple part and
  $\nsf()$ refers to the number of simple factors of a
  semisimple group.
  Then there exists a positive-dimensional normal subgroup $N$ of $G$ such that:
  \begin{enumerate}
  \item The fixed point set $M^N$ is non-empty
    (and $G$-invariant); let $B$ be a component of~$M^N$ containing
    principal $G$-orbits.
  \item $B/G$ has empty boundary and is contained in all faces of $X$.
  \item In particular:
    \begin{enumerate}
    \item[(i)] $N$ contains, up to conjugation, all isotropy groups of $G$
    corresponding to orbit types of strata of codimension one in $X$.
  \item[(ii)]
At a generic point of $B$, the slice representation of $N$ has orbit
    space with non-empty boundary.
  \item[(iii)] If, in addition, $M$ is simply-connected,
then the statement in (ii) is true with $N$ replaced by its
identity component $N^0$.
    \end{enumerate}
  \end{enumerate}
\end{thm}

This theorem will be proved in section~\ref{pf-main}. A rather
straightforward
modification of the argument proves a strengthened version
in which $M$ is only assumed to have positive $k$-th Ricci curvature,
inequality~(\ref{ineq}) is assumed to hold with $n$ replaced by $n-k+1$ and
the same conclusions are derived. Recall that a Riemannian manifold $M$ has
\emph{positive $k$-th Ricci curvature} if for each $p\in M$ and any $k+1$
orthonormal tangent vectors $e_0$, $e_1,\ldots,e_k$ at~$p$, the sum of
sectional curvatures $\sum_{i=1}^k K(e_0,e_i)>0$~\cite{Wu}.
The main examples with $k>1$ are compact
locally symmetric spaces with rank $\geq2$.

The following corollary of Theorem~\ref{main} is an immediate
consequence of~\cite[Theorem~7]{Wi2}.

\begin{cor}
  The orbit space $X$ is homeomorphic to the join of an $(f-1)$-simplex and the
  space (containing $B$) given by the intersection of all faces, where
  $f\leq\dim X$ is the number of faces of $X$.
\end{cor}

\subsection{Outline of proof of Theorem~\ref{main}}
The basic idea is to construct a certain normal subgroup of $G$
that contains all isotropy groups
associated to codimension one strata of~$X$ and prove that its
fixed point set is non-empty.
Suppose first $G$ is a simple Lie group. An involutive inner
automorphism of $G$ defines a symmetric space of inner type $G/K$ and
indeed corresponds to the geodesic symmetry at the base point of $G/K$.
On one hand, we can estimate the codimension of the fixed point set
of the involution in $M$, if we choose it to fix a regular point or an
important point (i.e.~a point projecting to a codimension one stratum
of $X$), which we can always do.
On the other hand, a finite number (which can be estimated
in terms of the geometry of $G/K$) of conjugates of the
involution generate a dense subgroup of $G$ (this is because
they correspond to geodesic symmetries of $G/K$ at generic points,
and these will generate sufficient transvections of $G/K$).
Combining these two observations yields, via Frankel's Theorem,
an estimate on the codimension of the fixed point set of $G$,
which is thus non-empty if the dimension of $M$ is sufficiently high.
In the case of a general compact connected Lie group, the
argument is more technical and one
proceeds by induction using the simple factors and the center.

\subsection{The Abelian case}
We illustrate some ideas in the proof in the much
simpler case of a torus action. So let a torus $T^k$
act effectively and isometrically on an orientable
connected complete $n$-manifold $M$ of positive sectional curvature
and assume $n\geq 2k$.
Note that the principal isotropy group $T_{pr}$ is trivial, since
it is a normal subgroup. If $p$ is an important point, $T_p$ is an
Abelian group that acts
simply transitively on the unit sphere of the non-trivial component
of the slice representation, and hence $T_p=S^0$ or $T_p=S^1$; the
first case cannot occur, as the non-trivial element in $T_p=S^0$
would act as a reflection on a codimension one hypersurface of $M$
and this is forbidden by the orientability of $M$.
We choose a point for each codimension one stratum in $X$
and end up with points $p_1,\ldots, p_\ell$. Let $L=T_{p_1}\cdots T_{p_\ell}$
be the group generated by the $T_{p_i}=S^1$.
Since $T$ is Abelian, the codimension of the fixed point set of $T_{p_i}$
is $2$. Owing to Frankel's Theorem, $\dim M^L\geq \dim M-2\dim L\geq2\dim T/L\geq0$,
so $M^L\neq\varnothing$. Let $\tilde B$ be a component of $M^L$
of maximal dimension. Now $T/L$ acts on $\tilde B$ and
$\dim\tilde B\geq2\dim T/L$.
If $\partial(\tilde B/T)\neq\varnothing$, we can repeat the
procedure; since $\dim T/L<\dim T$, the procedure must eventually stop.
We obtain a subtorus $S$ of $T$ containing $L$ and hence
all isotropy groups of codimension one strata of $X$,
whose fixed point
set $M^S$ has a component $B$ such that $\partial(B/T)=\varnothing$.

\subsection{Example}
Let $T^2=S^1_1\times S^1_2$ act on $M=S^5(1)$ by
$(S^1_1,\R^2)\times(S^1_2,\R^2\oplus\R^2)$, namely,
(standard action)$\times$(Hopf action). Then $X=S^3_+(\frac12)$,
$\partial X=S^2(\frac12)$, $N=S_1^1$, $B=M^N=S^3(1)$ and
$B/T^2=\partial X$.

\subsection{Structure of the paper}
After a short section on preliminaries, we show in section~\ref{sec:nice}
that the presence of boundary in the orbit space of the action
implies the existence of certain \emph{nice involutions},
whose codimension of the fixed point set we can estimate (Lemma~\ref{nice}),
unless some special situation occurs.
This is followed by section~\ref{tg} in which a
problem of independent interest
about the geometry of symmetric spaces is investigated,
namely, we want to know how many geodesic symmetries
of a compact symmetric space are needed to generate a
dense subgroup of the transvection group (compare
Proposition~\ref{dense} and Theorem~\ref{ell}).
In sections~\ref{pf-simple} and~\ref{pf-main}, we
apply the results of the two previous sections to
prove Theorems~\ref{simple} and~\ref{main}, respectively.
Section~\ref{reducible}
is devoted to establishing conditions under which a reducible
representation can have orbit space with non-empty
boundary (Proposition~\ref{prop:reducible} and Corollary~\ref{cor:reducible}).
The proofs of our applications are
finally collected in section~\ref{appl}.

The authors wish to thank Alexander Lytchak, Ricardo Mendes
and David Gonz\'alez-\'Alvaro for fruitful discussions and
valuable comments, and the anonymous referee for
constructive comments and recommendations which helped 
us to significantly improve the presentation.
Part of this work was completed while the first author was visiting
  the University of Cologne; he would like to thank
  Alexander Lytchak for his hospitality.

\section{Preliminaries}\label{prelim}

Let $G$ be a compact Lie group of isometries of a connected complete
orientable Riemannian manifold $M$.
Let $X$ be the orbit space $M/G$ equipped with the
induced quotient metric. We generally assume that the action is effective.

The subset of  $M$ consisting of all points with isotropy groups
conjugate to $G_p$ is a submanifold of $M$, denoted by $M_{(G_p)}$, called
an \emph{isotropy stratum} of $M$, and
projects to a Riemannian totally geodesic submanifold of $X$
denoted $X_{(G_p)}$,
called an \emph{isotropy stratum} of $X$, which contains
the point $x=Gp$.

Locally at $p\in M$, the orbit decomposition of $M$ is completely
determined by the \emph{slice representation} of $G_p$ on the normal space
$\nu_p(Gp)$. The set of $G_p$-fixed vectors
in $\nu_p(Gp)$ is tangent to $M_{(G_p)}$, and the action on its
orthogonal complement in $\nu_p(Gp)$ has
cohomogeneity equal to the codimension of $X_{(G_p)}$ in $X$.

A point $p\in M$ is called \emph{regular} if the slice representation
at $p$ is trivial. It is called \emph{exceptional} if it is not regular
and the slice representation has discrete orbits. If it is neither regular
nor exceptional, it is called \emph{singular}. The set $M_{reg}$ of all
regular points in $M$ is open and dense, and $X_{reg}$ is connected and
convex. $X_{reg}$ is the stratum corresponding to the unique
conjugacy class of minimal appearing isotropy groups; these are called
\emph{principal isotropy groups}.

The \emph{boundary} of $X$ is the closure of the union of all
strata of codimension $1$ in $X$. It is denoted by $\partial X$.
A point $p\in M$ is projected to a stratum of codimension $1$ in $X$
if and only if the non-trivial component of the slice representation
has cohomogeneity $1$; we will call such points \emph{$G$-important}.

We recall the easy but perhaps not much noticed fact that
the components of the fixed point set of a connected group of isometries of
an orientable manifold are \emph{orientable}
(closed totally geodesic) submanifolds~\cite[Theorem~3.5.2]{Z}.

\section{Nice involutions}\label{sec:nice}

Under the assumptions of section~\ref{prelim},
a \emph{nice involution} is a non-central element $\sigma\in G$ whose square
$\sigma^2$ is in the center of $G$ and whose fixed point in $M$ is
non-empty and has a component of
codimension at most $c+\dim G/K$, with $c\leq4$,
where $K=G^\sigma$ is the centralizer
of $\sigma$.
Nice involutions will play an important role
in estimating the codimensions of fixed point sets of certain groups
of isometries of~$M$.

Regarding the terminology in the statement of the next
lemma, recall that, in the case of finite principal isotropy groups, 
along each component of a 
codimension one stratum of the orbit space,
the connected slice representation is equivalent to
one of $(\Z_2,\R)$, $(S^1,\C)$ or $(S^3,\Q)$, up 
to a trivial subrepresentation (see the proof 
of the lemma and compare~\cite[section~4]{GL}); 
we call the corresponding 
component respectively of $\Z_2$-, $S^1$- or $S^3$-type. 

\begin{lem}\label{nice}
Assume $G$ is a compact connected Lie group
and $\partial X\neq\varnothing$.
Then nice involutions exist unless the following
situation $(\mathcal S)$ is present:
the principal isotropy group is finite of odd order, all
boundary components of $X$ are of $S^1$-type, and the
identity components of their
isotropy groups are contained in the center of~$G$.
\end{lem}

\Pf Assume the situation $(\mathcal S)$ does not happen.
We will look for certain 
involutions $\sigma\in G_p$ such that
$p\in M$ projects to a stratum
of codimension at most $1$ in~$X$, and later 
estimate the codimension of their fixed point sets, proving that 
they are nice involutions. As we will see, in certain cases
there are different kinds of possible choices for~$\sigma$. 

Fix a principal isotropy group $H$. 
Note that any element of $H$
which is central in~$G$ belongs to all principal isotropy groups
and thus lies in the kernel of the $G$-action, which
we have assumed to be trivial.

Assume first that 
$H$ is finite. For any $G$-important point $p\in M$, the
isotropy group $G_p$ acts transitively on the unit sphere $S^a$
in the non-trivial component of the slice representation.
It follows that $G_p/H$ is diffeomorphic to $S^a$.
In particular, for $a\geq1$ there is a finite covering $G_p^0\to S^a$.
If $a\geq2$, $S^a$ is simply-connected, so the covering is a
diffeomorphism and thus $a$ equals $3$. We can take $\sigma=-1\in S^3\approx G_p^0$,
or a square root of this element if it is central
in~$G$. If $a=1$, then $G_p^0$ is a finite
covering of $S^1$, hence, $G_p^0\approx S^1$ and we may assume this 
subgroup is non-central (thanks to our assumption that 
$(\mathcal S)$ is not present). 
Then $Z(G)\cap G_p^0$ is at most a cyclic group
and again we can take~$\sigma$ to be a square root of
an element of $Z(G)\cap G_p^0$.
If $p\in M$ is a
$G$-important point with $G_p/H\approx S^0=\Z_2$, there is $\sigma'\in G_p^0$
acting as~$-1$ on the $1$-dimensional non-trivial component of the
slice representation. Since~$G$ is connected and $M$ is orientable,
$\sigma'$ cannot be central. Since $(\sigma')^2$ is
trivial on the slice, it is an element of $H$.  
In case $H$ has odd order, also $(\sigma')^2$ has odd order,
say $2b+1$, and we can take $\sigma=(\sigma')^{2b+1}$.

If $H$ is finite with even order or $\dim H>0$, 
it is clear that
we can find an element $\sigma\in H$ of order $2$.
If possible, this is our preferred choice of~$\sigma$,
in terms of obtaining a better estimate on the 
codimension of $M^\sigma$, see below. 
 
It remains only to estimate the
codimension of the fixed point set of~$\sigma$. The non-trivial component of the
slice representation at $p$ has dimension $c$ equal
to~$0$, $1$, $2$ or $4$ according to whether $p$ is a regular point
or an important point projecting to a boundary stratum of type
$\Z_2$, $S^1$ or $S^3$, respectively. Along the tangent space $T_p(Gp)$,
the codimension of the fixed point set of $\sigma$ is bounded
by the dimension of the $(-1)$-eigenspace of $\mathrm{Ad}_\sigma$,
that is, $\dim G/K$. Hence the component through~$p$ of the fixed point
set $M^\sigma$ has codimension at most $c+\dim G/K$. \EPf

\begin{rmk}\label{rmk:nice}
The proof of the lemma shows that
we can take $c=0$ if $\dim H>0$ or $H$ is finite with even order,
$c=2$ if $H$ is finite and there are
no $S^3$-boundary components, and $c=4$ in general. 
\end{rmk}

\section{Generic totally geodesics submanifolds of
  compact symmetric spaces}\label{tg}

The second tool that we will use is an invariant attached to
a compact connected symmetric space~$M$.
Define $\ell_M$ to be the
minimal number $\ell$ such that there exists $p_1\dots,p_\ell\in M$
``spanning'' $M$, in the sense that these points do not lie
in a proper connected closed totally geodesic submanifold
of~$M$.

More specifically, let $p_1,\ldots,p_k\in M$ with $k\geq2$
be generic points in the sense that each pair $(p_i,p_j)$ with $i\neq j$
is connected by a unique shortest geodesic. In this case, it is
clear that the intersection of all closed connected totally
geodesic submanifolds of $M$ which contain~$p_1,\ldots,p_k$
has a connected component
containing $p_1,\ldots,p_k$, which we call the
\emph{span} of $p_1,\ldots,p_k$ and denote by
$\langle p_1,\ldots,p_k\rangle$. It is easy to see
that $\ell_M$ is the minimal number $\ell$ such that
there exist generic points $p_1,\ldots,p_\ell\in M$
with $\langle p_1,\ldots,p_\ell\rangle=M$.

Note also that $\langle p_1,\ldots,p_k\rangle$ equals
$\exp_{p_1}(T_{p_1}\langle p_1,\ldots,p_k\rangle)$,
where $T_{p_1}\langle p_1,\ldots,p_k\rangle$
coincides with the
intersection of all Lie triple systems in $T_{p_1}M$ containing
$v_2,\ldots,v_k$, where $v_i$ is tangent to the geodesic joining
$p_1$ to $p_i$. We deduce that $\ell_M=\ell_{M'}$ for a
Riemannian covering $M'\to M$.

%Next, we impose further generic conditions on the points.
Consider now the case $k=2$. It is clear that
$\langle p_1,p_2\rangle$ is either a closed geodesic
or the closure of a non-periodic infinite geodesic, that is,
in any a case a flat torus
of~$M$. The extreme case occurs when
$p_2$ is a regular point with respect
to the isotropy action at $p_1$, and the
geodesic through $p_1$ and $p_2$ is dense in the
unique maximal flat torus $T_{12}$ of $M$ containing those
points; in this case $\langle p_1,p_2\rangle=T_{12}$.
In particular $\langle p_1,\ldots,p_k\rangle$ for $k\geq2$ has maximal rank
and $\ell_M\geq3$ if $M$ is not flat.
Indeed we shall see that $\ell_M=3$ if $M$ is an irreducible
compact symmetric space of inner type, unless $M$ is one of
$\Q P^2$, $\Ca P^2$ %, $\E6/\F$
or $\mathsf{Gr}_k(\K^n)$ with $n>3k$
(here $\K=\R$, $\C$ or $\Q$).

From now on, we impose further genericity conditions. Let
 $p_1,\ldots,p_k\in M$ with $k\geq2$
be generic points in the sense that each pair $(p_i,p_j)$ with $i\neq j$
is connected by a unique shortest geodesic, $p_j$ is regular with
respect to the isotropy action at $p_i$ so that $p_i$ and $p_j$ are contained
in a unique maximal flat torus $T_{ij}$ of $M$, and the product
of the geodesic symmetries of $M$ at $p_i$, $p_j$ is a transvection
generating a group acting transitively on a dense subset of $T_{ij}$.
With these genericity conditions,
denote by $L=L_{p_1,\ldots,p_k}$ the
closure of the group consisting of even products of
geodesic symmetries of $M$ at $p_1,\ldots,p_k$.
Then $L$ is connected and $L(p_1)=\cdots=L(p_k)$ is a submanifold
of $\langle p_1,\ldots,p_k\rangle$. Since the geodesic
symmetry of $M$ at any point of $L(p_1)$ leaves this
submanifold invariant, $L(p_1)$ is totally geodesic and
hence $\langle p_1,\ldots,p_k\rangle=L(p_1)$.

Write $M=G/K$ where $G$ is the identity component of the
isometry group of~$M$ (the transvection group of $M$)
and $K=G_{p_1}$, and let $\Lg=\Lk+\Lp$
be the decomposition into the $\pm1$-eigenspaces of the involution
induced by the geodesic symmetry at $p_1$;
here $\Lk$ is the Lie algebra of $K$ and $\Lp\cong T_{p_1}M$.
Fix a maximal Abelian
subspace $\La$ of $\Lp$. Then there are decompositions
\begin{equation}\label{root-space-decomp}
  \Lk=\Lk_0+\sum_{\lambda\in\Lambda^+}\Lk_\lambda,
  \quad\Lp=\La+\sum_{\lambda\in\Lambda^+}\Lp_\lambda,
  \end{equation}
where $\Lambda$ is a (possibly non-reduced) root system. The
\emph{marked Dynkin diagram of $M$} is the
Dynkin diagram of $\Lambda$ where each vertex is labeled
by the multiplicity~$m_\lambda=\dim\Lp_\lambda$, with a special
rule in case of non-reduced roots, see~\cite[p.~118]{Loos2}.

Suppose $p_1$, $p_2\in M$ are generic points and
$\langle p_1,p_2\rangle=\exp_{p_1}(\La)$. A generic
choice of $p_3\in M$ corresponds under $\exp_{p_1}$
to~$v_3\in T_{p_1}M\cong\Lp$ such that all $\Lp_{\lambda}$-components
of~$v_3$ are
nonzero, since the linear isotropy representation of $K$ on $\Lp$ preserves
the $\Lp_\lambda$.
It follows that for generic $p_1,\ldots,p_k\in M$
the marked Dynkin diagram of
$\langle p_1,\ldots,p_k\rangle$ as a symmetric space already
coincides with that of $M$ if $k=3$, up to the multiplicities which are bounded
above by those of $M$, and hence also for $3< k<\ell_M$.

Below we compute $\ell_M$ for compact irreducible symmetric spaces
of inner type. The reducible case is covered
by the following lemma.

\begin{lem}\label{product}
  Let $M_1$ and $M_2$ be two compact symmetric spaces.
  Then
  \[ \ell_{M_1\times M_2}=\max\{\ell_{M_1},\ell_{M_2}\}. \]
  \end{lem}

\Pf Given generic points $p_1,\ldots,p_k\in M_1\times M_2$
with $k\geq2$, the closed totally geodesic submanifold
$\langle p_1,\ldots p_k\rangle$ has maximal rank in $M_1\times M_2$,
so it is of the form $N_1\times N_2$, where $N_i$ is a totally
geodesic submanifold of $M_i$. The result follows. \EPf

\medskip

The importance of the invariant $\ell_M$ for us lies in the
following proposition.

\begin{prop}\label{dense}
Let $G$ be a compact connected Lie group and assume $\sigma\in G$
is a non-central
element whose square is central. Let $K=G^\sigma$ be the centralizer of $\sigma$.
Assume $G$ acts almost effectively on the symmetric space $M=G/K$.
Then there exist $g_1,\ldots,g_{\ell_M}$ such that
the group generated by $g_1\sigma g_1^{-1},\ldots, g_{\ell_M}\sigma g_{\ell_M}^{-1}$
is dense in~$G$.
\end{prop}

\Pf Note that the assumption that $G$ acts almost effectively on $M$
says that $\sigma$ does not centralize a normal subgroup of $G$
of positive dimension and $G$ is the identity component of the
isometry group of $M$, up to a finite covering.
Let $p$ denote the base point of $M=G/K$ and choose
$g_1,\ldots,g_k$ such that $g_1p,\ldots, g_kp\in M$ are in generic
position. Then the totally geodesic submanifold
\[ N:=\langle g_1p,\ldots,g_kp\rangle \]
is closed and connected, and the closure of the
group generated by $g_1\sigma g_1^{-1},\ldots, g_k\sigma g_k^{-1}$ is a
closed subgroup of $G$
containing all transvections of $N$;
%therefore its action on $N$
%coincides with the identity component of the isometry group of $N$.
But $N=M$ for $k=\ell_M$.\hfill\mbox{ } \EPf

\medskip

For the sake of computation of $\ell_M$, we next introduce another
invariant of a compact symmetric space $M=G/K$, where $G$ is the
transvection group of $M$. Define  $h_M$ to be the maximal number $h$ such that the
principal isotropy group~$H$ of
the diagonal action of $G$ on the $h$-fold product $M^h$ is
non-trivial.
Note that for $h=1$ the principal isotropy group is $K$, and for $h=2$
the principal isotropy group is the principal isotropy group $K_{pr}$
of the linear isotropy representation of $K$
on the tangent space $T_pM$ at the base point~$p$, which is never trivial,
so $h_M\geq2$. Indeed $h_M=1+\tilde h_M$, where $\tilde h_M\geq1$ is the
maximal number $\tilde h$ such that
the principal isotropy group of the diagonal action of $K$ on the $\tilde h$-fold
sum $\oplus^{\tilde h}T_pM$ is non-trivial.

\begin{prop}\label{h}
Let $M=G/K$ and $H$ be as above. Then $h_M+1\leq\ell_M\leq h_M+2$.
Furthermore, in case $\ell_M=h_M+2$ there is a closed
connected totally geodesic submanifold $N_2$ of $M$ (different from $M$)
of codimension at most $\dim H$ such that
$N_2=\langle p_1,\ldots,p_{\ell_M-1}\rangle$ for generic
points $p_1,\ldots,p_{\ell_M-1}\in M$. In particular, if
$H$ is finite then $\ell_M=h_M+1$.
\end{prop}

\Pf Given
generic points
$p_1,\ldots,p_{h}\in M$, $h=h_M$, they all lie in $M^H$, up to
replacing $H$ by a conjugate group. Note that $M^H$ is a closed
totally geodesic submanifold of $M$, but it is not
necessarily connected. However, $H$ centralizes the
geodesic symmetries at the $p_i$ and hence centralizes
the group $L_{p_1,\ldots,p_{h}}$. It follows that
$\langle p_1,\ldots,p_h\rangle\subset M^H$.
Since the former submanifold is connected, we deduce that
$\ell_M-1\geq h_M$.

It remains to obtain the upper bound for $\ell_M$. Assume $\ell=\ell_M>h+1$ and
fix generic points $p_1,\ldots,p_{\ell-2},q_{\ell-1}\in M$.
We have the
closed connected totally geodesic
submanifolds $N_1=\langle p_1,\ldots,p_{\ell-2}\rangle$
and $N_2=\langle p_1,\ldots,p_{\ell-2},q_{\ell-1}\rangle$.
Owing to the fact that
the number of closed connected totally geodesic submanifolds
of a compact symmetric space, up to congruence, is countable,
there is a subset of
positive measure $U$ of $M$ such that
for all $p_{\ell-1}\in U$, the flag of closed connected
totally geodesic submanifolds
$N_1\subset\langle p_1,\ldots,p_{\ell-1}\rangle$ is $G$-conjugate
to $N_1\subset N_2$. In other words, for all  $p_{\ell-1}\in U$
there is $\iota\in G$ such that $\iota(p_i)\in N_1$ for
$i\leq\ell-2$ and $\iota(p_{\ell-1})\in N_2$.

The isometry group $\mathrm{Iso}(N_1)$ is a compact Lie group with finitely
many connected components. By~\cite[Lemma~7.5]{Wi}, there is a
finite subgroup~$F$ of $\mathrm{Iso}(N_1)$ meeting every component.
We can find 
$\psi$ in the identity component $\mathrm{Iso}_0(N_1)$
such that $\psi\cdot\iota|_{N_1}\in F$. Every geodesic symmetry of $N_1$
uniquely extends to a geodesic symmetry of $N_2$
and then to a geodesic symmetry of $M$,
and hence every transvection of $N_1$ admits an extension (not necessarily
unique) to a transvection of $N_2$ and then to a transvection of $M$.
%There are natural
%inclusions of transvection groups of symmetric spaces
%$\mathrm{Iso}_0(N_1)\subset\mathrm{Iso}_0(N_2)\subset G$, so we
Since the group generated by transvections at a fixed point 
of a symmetric space coincides with the identity component 
of the isometry group of the symmetric space, 
we may consider $\psi\in G$ and then the element $\psi\cdot\iota\in G$
also maps $p_{\ell-1}$ to~$N_2$. We have shown that we can always take
$\iota\in\hat H$, where
\[ \hat H:=\{g\in G\;|\;\mbox{$g(N_1)=N_1$ and $g|_{N_1}\in F$}\}. \]

Notice that $\hat H$ is closed subgroup of $G$
and hence a Lie group, with the same identity
component as the isotropy group $\bar H$ of $(p_1,\ldots,p_{\ell-2})$.
Since $U$ has positive measure and $\hat H(N_2)\supset U$, it follows
that the codimension of $N_2$ is bounded above by $\dim\hat H=\dim \bar H$.
In particular $\bar H$ is nontrivial and hence $h=\ell-2$ and
$\bar H=H$. \hfill\mbox{ }\EPf

\begin{cor}\label{l}
If every maximal connected closed totally geodesic submanifold
of $M$ is given as a component of the fixed point set of a
subgroup of $G$, then $\ell_M=h_M+1$.
\end{cor}

\Pf Given generic points $p_1,\ldots,p_{\ell_M-1}\in M$,
there exists a connected closed totally geodesic submanifold
$N$ containing those points, and we can assume
$N$ is maximal. By assumption, $N$ is a component of the fixed
point set of a non-trivial subgroup $H$ of $G$.
Now $H$ is contained in the isotropy group $\tilde H$
of a generic $(\ell_M-1)$-tuple of points of $M$, so
$h_M\geq\ell_M-1$. \EPf

\begin{rmk}
  A connected totally geodesic submanifold of $M$ is called
\emph{reflective} if it is a connected component of the fixed
point set of an involutive isometry of $M$; if, in addition,
the involutive isometry can be taken in the transvection
group of $M$, then the submanifold will be called \emph{inner reflective}.
It follows from Corollary~\ref{l} that
if every maximal connected closed totally geodesic submanifold
of $M$ is inner reflective, then $\ell_M=h_M+1$.
\end{rmk}

\begin{thm}\label{ell}
  The invariant $\ell_M$ for various irreducible symmetric
  spaces $M$ of compact type is listed in Table~3, including all
  spaces of inner type.
\end{thm}

\begin{table}[h]
  \[ \begin{array}{|c|c|c|}
    \hline
   M=G/K & \ell_M & \text{Conditions} \\
   \hline
   \SO n/(\SO p \times \SO{n-p}) & & p\leq n/2 \\
   \SU n/\sf S(\U p \times \U{n-p}) & \max\{3,\lceil\frac np\rceil\}& p\leq n/2 \\
   \SP n/(\SP p \times \SP{n-p}) & & p\leq n/2,\ (n,p)\neq(3,1) \\
   \hline
   \SP 3/\SP 1 \times \SP{2} &\multirow{2}{*}{\rm 4} &\multirow{2}{*}{$-$} \\
   \F/\Spin9&& \\
   \hline
   \SP n/\U n &&-\\
   \SO{2n}/\U n&&n\geq5\\
   \G/\SO4&&-\\
   \F/\SP3\SP1&&-\\
   \E6/\Spin{10}\U1   &&-\\
   \E6/\SU6\SU2&3&-\\
   \E7/\E6\U1&&-\\
   \E7/(\SU8/\Z_2)&&-\\
   \E7/\Spin{12}\SU2&&-\\
   \E8/\Spin{16}&&-\\
   \E8/\E7\SU2&&-\\
   \hline
  \end{array}\]
  \smallskip
  \begin{center}
    \sc Table~3: The invariant $\ell_M$ for some irreducible symmetric
    spaces of compact type.
  \end{center}
\end{table}

\Pf We run through the cases.
\iffalse
It is useful to recall the \emph{index} of
a symmetric space $M$, introduced by Onishchik~\cite{O}
as the minimal codimension of a proper closed
totally geodesic submanifold of~$M$.
More recently, Berndt and Olmos~\cite{BO1,BO2} have
extended his results to a larger group of symmetric spaces.
\fi

\emph{Symmetric spaces of maximal rank.}
They are
\begin{gather*}
\SO{2p}/(\SO p \times \SO p),\ \SU{n}/\SO n,\ \SP n/\U n,\\
\SO{2p+1}/(\SO{p+1}\times\SO p),\ \E6/(\SP4/\Z_2),\
\E7/(\SU8/\Z_2),\\ \E8/\mathsf{SO}'(16),\ \F/\SP3\SP1\ \mbox{and}\ \G/\SO4
\end{gather*}
(not all listed in Table~3).
The condition $\mathrm{rk}\,M=\mathrm{rk}\,G$ is equivalent
to the effective $K_{pr}$ being finite~\cite[Proposition~4.1]{Loos2} (and
indeed isomorphic to $\Z_2^{\mathrm{rk}\,M}$).
Therefore $h_M=2$ and $\ell_M=3$ by Proposition~\ref{h}.

\iffalse
\emph{The quaternionic projective plane~$\Q P^2=\SP3/(\SP1\times\SP2)$.}
We have $h_M=3$, and then $\ell_M=4$ by Corollary~\ref{l}.
In fact a maximal closed connected
totally geodesic
submanifold of $\Q P^2=\SP3/(\SP2\times\SP1)$ is either
$\C P^2=\U3/(\U1\times\U2)$
or $\Q P^1=(\SP2\times\SP1)/(\SP1\times\SP1\times\SP1)$~\cite{Wo2} and
hence inner reflective, as
$\U3$ and $\SP2\times\SP1$ are components of fixed point sets of an inner
automorphism of $\SP3$.
One can also directly check
that any three points in $\Q P^2$ are contained in a
totally geodesic $\C P^2$.
\fi

\emph{The Cayley projective plane~$\Ca P^2=\F/\Spin9$.}
The linear isotropy representation of $K=\Spin9$ on $\R^{16}$
has $K_{pr}=\Spin7$ and  corresponding $K_{pr}$-irreducible
decomposition
$\R\oplus\R^7\oplus\R^8$. The principal isotropy group of this action is $H\cong\SU3$,
with corresponding decomposition $4\R\oplus2\C^3$. The principal isotropy group
of this action is trivial, so $h_M=3$.
A maximal closed connected totally geodesic
submanifold of $\Ca P^2=\F/\Spin9$ is either
$\Ca P^1=\Spin9/\Spin8$ or~$\Q P^2=(\SP3\cdot\SP1)/(\SP2\cdot\SP1\cdot\SP1)$.
Since $\Spin9$ and $\SP3\cdot\SP1$ are components of fixed point sets of
inner automorphisms of $\F$, $\Ca P^1$ and $\Q P^2$ are inner
reflective. We deduce from Corollary~\ref{l}
that $\ell_M=4$.

\emph{Grassmann manifolds.} Let $M=\mathsf{Gr}_r(\K^n)$ with $n\geq2r$,
where $\K=\R$,
$\C$ or $\Q$.
Given $p_1,\ldots,p_k\in M$,
these points respectively
lift to $r$-dimensional $\K$-subspaces $\pi_1,\ldots,\pi_k\subset \R^n$.
If $k<\frac nr$, then clearly the span of $\pi_1,\ldots,\pi_k$ is
a proper subspace of $\K^n$ and $p_1,\ldots,p_k\in\mathrm{Gr}_q(\K^{kr})$,
so that $\ell_M-1\geq k$. Note that we can always
take $k=m-1$, where $m:=\lceil \frac nr \rceil\geq2$, so $\ell_M\geq\max\{3,m\}$.

If $m=2$, then $n=2r$ (in the case
$\K=\R$, $M$ is a space of maximal rank and this case has already been
examined) and
it is not so difficult to find three points in $M$ not contained
in a proper connected closed totally geodesic submanifold,
implying $\ell_M=3$.
Next we assume $m\geq3$ and want to show that there exist
$p_1,\ldots,p_m\in M$ such that $N:=\langle p_1,\ldots,p_m\rangle$
coincides with~$M$. This will prove $\ell_M=m$.
Let $\{e_i\}_{i=1}^n$ be the canonical $\K$-basis of
$\K^n$ and
consider $p_1,\ldots,p_m$ associated to the $r$-dimensional subspaces
(note that $(m-1)r+1\geq n-r+1$):
\begin{eqnarray*}
  \pi_1&=&\mathrm{span}(e_1,\ldots,e_r),\\
  &\vdots&\\
\pi_{m-1}&=&\mathrm{span}(e_{(m-2)r+1},\ldots,e_{(m-1)r}),\\
    \pi_m&=&\mathrm{span}(e_{n-r+1},\ldots,e_n).\\
\end{eqnarray*}
By slightly perturbing the points $p_1,\ldots,p_m$, we can ensure
that $N$ is a connected closed
totally geodesic submanifold of maximal rank and same restricted
root system as~$M$. In cases $\K=\R$ or $\C$,
using the classification~\cite[p.~119 and p.~146]{Loos2}
this already implies that
$N$ is a $\K$-Grassmannian. In case $\K=\Q$, below
we distinguish between $r>1$ and $r=1$ to prove
that $N$ is an $\Q$-Grassmannian.
In any case, since $\{e_i\}_{i=1}^n$ has been perturbed to another
$\K$-basis of $\K^n$, we must have~$N=M$.

In case $\K=\Q$ and $r>1$ we check that, for generic $p_1,\ldots, p_m$,
$N$ is an $\Q$-Grassmannian as follows. Consider the restricted
root space decompostion~(\ref{root-space-decomp}) where
$p_1$ is the basepoint and $\langle p_1,p_2\rangle=\exp_{p_1}(\La)$.
It is not difficult to see that $v_3\in T_{p_1}M\cong\Lp$ can be chosen so that
$\langle p_1,p_2,p_3\rangle$ is not a $\C$-Grassmannian.
Note also that $N=\SO{4r+2}/\U{2r+1}$ is not a totally geodesic
submanifold of $M$ (one way to see that is as follows:
consider the restricted root system
$\{\theta_i\pm\theta_j,\theta_i,2\theta_i\}$ of type $\mathsf{BC}_r$;
of course
$[\Lp_{\theta_1+\theta_2},\Lp_{\theta_1-\theta_2}]\subset \Lk_{2\theta_1}+\Lk_{2\theta_2}$;
  one computes directly that the left hand-side has
  dimension $3$ in case of $M$, and hence also in case of $N$
  as the multiplicities of $\theta_1\pm\theta_2$ equal~$4$ in both cases;
  however the right hand-side has dimension $2$ in case of $N$).
\iffalse
(in case $r=2$ this follows from~\cite[Table~VIII]{CN}
or~\cite[Table~4]{BO2}, and in case $r\geq3$ this follows
from~\cite[Proposition~7.4]{BO2} and~\cite[p.~173]{Le3}).
\fi
It follows from  the
classification~\cite[p.~119 and p.~146]{Loos2} that
$\langle p_1,p_2,p_3\rangle$ is an $\Q$-Grassmannian
and so is $\langle p_1,\ldots,p_m\rangle$.

In case $\K=\Q$ and $r=1$, $M=\Q P^{n-1}$ is of type ${\sf BC}_1$ and $m=n$.
It is not difficult to see that for a generic choice of points,
$\langle p_1,p_2,p_3\rangle=\C P^2$. If $n=3$,
this is a maximal totally geodesic
submanifold, so $\ell_{\mathbb HP^2}=4$. If $n>3$,
$\langle p_1,p_2,p_3,p_4\rangle$ is an~$\Q$-projective space,
and so is  $\langle p_1,\ldots,p_m\rangle$. We finish as above
to deduce that $\ell_M=m=n$.

\iffalse
If $s=1$ then $N_1$ is maximal in $M$,
so $\ell_M=q$. If $s\geq2$, then $N_1$ is the fixed point set of $H$,
the lower $\SO s$-block of $\SO n$ and, according to the
argument in the proof of Proposition~\ref{h}, $\hat H(N_2)$ contains a
subset of positive measure of $M$,
\fi

\emph{The space $\SO{4n}/\U{2n}$.}
Note that the cases $n=1$ and $n=2$ are respectively locally
isometric to a sphere and a real Grassmannian, so we may assume $n\geq3$.
The linear isotropy representation is $\Lambda^2\C^{2n}$ with
$K_{pr}=\SU2^n$, so $h_M=2$. If $\ell_M=4$ then due to Proposition~\ref{h}
$M$ contains
a connected closed totally geodesic submanifold $N_2$ with
the same Dynkin diagram, and dimension at least
$(4n^2-2n)-3n=4n^2-5n$. According to~\cite[p.~119 and~146]{Loos2},
the submanifold with same diagram, not larger multiplicities
and maximal dimension
is $\mathsf{Gr}_n(\C^{2n})$, which has dimension $2n^2<4n^2-5n$ for $n\geq3$.
Hence $\ell_M=3$.

\emph{The space $\SO{4n+2}/\U{2n+1}$.}
Note that the case $n=1$ is locally
isometric to a $\C P^3$, so we may assume $n\geq2$.
The linear isotropy representation is $\Lambda^2\C^{2n+1}$ with
$K_{pr}=\SU2^n\U1$, so $h_M=2$. If $\ell_M=4$ then due to Proposition~\ref{h}
$M$ contains
a connected closed totally geodesic submanifold $N_2$
the same Dynkin diagram, and dimension at least
$(4n^2+2n)-(3n+1)=4n^2-n-1$. Note that $\SU{2n+2}$ is not a subgroup of
$\SO{4n+2}$ so, according to~\cite[p.~119 and~146]{Loos2},
the only candidate with same diagram, not larger multiplicities
and maximal dimension is $\mathsf{Gr}_n(\C^{2n+1})$
which, however, has dimension $2n^2+2n<4n^2-n-1$ for $n\geq2$.
Hence~$\ell_M=3$.

\emph{The space $\E6/\Spin{10}\U1$.}\label{evi}
The linear isotropy representation is $\C^{16}\otimes_{\mathbb C}\C$
with $H=K_{pr}=\U4$ and corresponding decomposition
$4\R\oplus2\C^4\oplus2\R^6$, so $h_M=2$.
If $\ell_M=4$ then, due to Proposition~\ref{h}, $M$ contains
a connected closed totally geodesic submanifold $N_2=\langle p_1,p_2,p_3\rangle$
of rank $2$, same Dynkin diagram, and dimension at least
$\dim M-\dim H=32-16=16$; here $p_1$, $p_2$, $p_3\in M$ are
generic points.
According to Chen and Nagano (cf.~\cite{CN,klein}), the
submanifolds under these conditions are $\mathsf{Gr}_2(\C^6)$
and $\SO{10}/\U5$.
The first submanifold is a connected component of the fixed point set
of the geodesic symmetry of $\E6/\SU6\SU2$ (compare~\cite[p.~1115]{klein}
and~\cite[Proposition~3.5]{K}).
Similarly, the second one is a polar submanifold, namely, a connected
component of the fixed point set of the geodesic symmetry of
$M$, see~\cite[p.~1119]{klein}. It follows that in both cases the isotropy
group in $\E6$ of a generic triple of points in $M$ is non-trivial,
which is a contradiction to $h_M=2$ (compare Corollary~\ref{l}).
Hence~$\ell_M=3$.

\emph{The space $\E6/\SU6\SU2$.}
The linear isotropy representation is $\Lambda^3\C^6\otimes_{\mathbb H}\C^2$
with $K_{pr}=T^2\cdot\Z_2$~\cite[p.~436]{HPT},
so $h_M=2$. If $\ell_M=4$, then due to Proposition~\ref{h}
$M$ contains
a connected closed totally geodesic submanifold $N_2$
of codimension~$2$, but this symmetric space is of
index bigger than~$2$~\cite[Thm.~1.2]{BO1}. Hence~$\ell_M=3$.

\emph{The space $\E7/\E6\U1$.}
The linear isotropy representation is $\C^{27}\otimes_{\mathbb C}\C$
with $K_{pr}=\Spin8$, and $h_M=2$. If $\ell_M=4$ then
$M$ contains a connected closed totally geodesic submanifold $N_2$
of $\sf C_3$-type, multiplicities bounded above by $(8,8,1)$
and dimension at least~$26$.
Looking at the list of diagrams~\cite[p.~119]{Loos2},
$N_2$ must be $\SO{12}/\U6$. This submanifold is a connected
component of the fixed point set of the involution of $M=\E7/\E6\U1$
induced by the involution of $\E7$ defining $\Spin{12}\SU2$ as
a symmetric subgroup (cf.~\cite[p.~70]{nagano} and~\cite[Proposition~3.5]{K});
since the latter is an
inner automorphism of $\E7$, it follows as in subsection~\ref{evi} that
the isotropy group of a generic triple of points of~$M$ is
non-trivial, a contradiction to~$h_M=2$.
Hence $\ell_M=3$.

\emph{The space $\E7/\Spin{12}\SU2$.}
The linear isotropy representation is $\C^{32}\otimes_{\mathbb H}\C^2$
with $K_{pr}=\Z_2^2\cdot\SP1^3$~\cite[p.~436]{HPT}, so $h_M=2$.
If $\ell_M=4$ then
$M$ contains a connected closed totally geodesic submanifold $N_2$
of rank $4$, $\F$-type, multiplicities bounded above by~$(4,4,1,1)$
and dimension at least $55$.
However, there exist no symmetric spaces under these
conditions~\cite[p.~119 and~146]{Loos2}.

\emph{The space $\E8/\E7\SU2$.}
The linear isotropy representation is $\C^{56}\otimes_{\mathbb H}\C^2$
with $K_{pr}=\Z_2^2\cdot\Spin8$~\cite[p.~436]{HPT}, so $h_M=2$.
If $\ell_M=4$ then
$M$ contains a connected closed totally geodesic submanifold $N_2$
of rank $4$, $\F$-type and dimension at least $84$.
Looking at the list of diagrams in~\cite{Loos2},
we see there are no submanifolds
under these conditions. Hence $\ell_M=3$. \EPf

\begin{rmk}\label{rmk:estimates}
For each symmetric space $G/K$ in Table~3,
a direct check shows that 
the following estimates hold: 
\begin{gather*}
 \ell_{G/K}\dim G/K\leq2\dim G;\\
\ell_{G/K}\leq 2\rk G+1. 
\end{gather*}
These will be used in the proof of Theorem~\ref{main}.
\end{rmk}

\section{Proof of Theorem~\ref{simple}}\label{pf-simple}

Although the proof of this theorem is contained
in the proof of Theorem~\ref{main}, to be proved in the 
next section, it is instructuve to do this proof first. 

Let $\mathcal L_G$ be given by~(\ref{lg}). 
Since $G$ is assumed to simple, situation~$(\mathcal S)$
in Lemma~\ref{nice} does not occur, so the lemma yields 
a nice involution $\sigma\in G$. Let $K=G^\sigma$. Then 
$G$ acts almost effectively on the symmetric space $G/K$. 
By Proposition~\ref{dense} we can find $g_1,\ldots,g_{\ell_{G/K}}\in G$, such that the group generated 
by $\sigma_i=g_i\sigma g_i^{-1}$ for $i=1,\ldots,\ell_{G/K}$ is dense 
in~$G$. Using Frankel's theorem and the codimension estimate 
for nice involutions, we obtain that 
\begin{eqnarray*}
\dim M-\dim M^G &=& \dim M-\dim M^{\sigma_1}\cap\cdots 
M^{\sigma_{\ell_{G/K}}} \\
&\leq& \sum_{i=1}^{\ell_{G/K}} \dim M - \dim M^{\sigma_i} \\
&\leq& \ell_{G/K}(4+\dim G/K) \\
&=& \mathcal L_G, 
\end{eqnarray*}
as desired. This completes the proof of the theorem.

\section{Proof of the main result}\label{pf-main}

We now proceed with the much more involved proof of
 Theorem~\ref{main}.
We follow a finite algorithm.
At each step, there are two possibilities, namely, the situation~$(\mathcal S)$
as in Lemma~\ref{nice} is present or not.

\subsection{$(\mathcal S)$ is not present}
By Lemma~\ref{nice}, we can choose a nice involution $\sigma\in G$.
Then $G/K$ is a symmetric space of inner type, where $K=G^\sigma$,
which locally splits as $G_1/K_1\times\cdots \times G_m/K_m$, where
each factor $G_i/K_i$ is not necessarily irreducible, but instead
$\ell_i:=\ell_{G_i/K_i}$ satisfy $\ell_i<\ell_{i+1}$ for
$i=1,\ldots m-1$. Furthermore we may take $G_i$ to be a
connected closed normal subgroup of $G$ acting with finite
kernel on~$G_i/K_i$.

Put $G'=G_1\cdots G_m\subset G_{ss}$, non-trivial connected
semisimple Lie group. Then $G=G'\cdot G''$ where
$G''$ denotes the identity component of the centralizer of $G'$ in $G$
and contains $Z(G)^0$. It follows that
$\alpha_G=\alpha_{G'}+\alpha_{G''}$ and 
$\beta_G=\beta_{G''}$. 

We will construct a 
component $\tilde B$ of $M^{G'}$ such that 
$\dim\tilde B>\alpha_{G''}+\beta_{G''}$.
The point here is, since $\tilde B$ is orientable and 
positively curved, in case 
$\tilde B/G''$ has non-empty boundary, we can check whether $(\mathcal S)$ 
is present or not and repeat the argument
for the action of $G''$ on $\tilde B$.
In repeating the argument for $G'$, we
get a similar decompostion $G'=(G')'\cdot(G')''$ and construct a component 
$\vardbtilde B$ 
of $(\tilde B)^{(G')'}$ on which $(G')''$ acts etc. Note that both $G'$ and $(G')'$ 
act trivially on $\vardbtilde B$. 
Since $\dim G''<\dim G$, this process will stop after finitely many steps. 

So each time that we repeat the argument,
we get a typical pair $(G',\tilde B)$. 
We construct an ascending chain of connected normal subgroups of $G$
by collecting the factors $G'$ in each step. The maximal 
element in this chain is a connected normal 
subgroup $G^\infty=G'\cdot(G')'\cdot((G')')'\cdots$ 
of $G$ that will contain all isotropy groups 
associated to codimension~$1$ strata of~$X$. On the other hand, 
as a fixed point set, in each step
the $\tilde B$ form a descending chain of 
connected totally geodesic submanifolds of $M$,
whose minimal element is a component $B$ of $M^{G^\infty}$,
which will be contained in all faces of $X$. Finally, we enlarge
$G^\infty$ to the (possibly disconnected) 
subgroup $N$ of $G$ that fixes $B$ pointwise.

Choose $\ell_m$ elements $g_1,\ldots,g_{\ell_m}\in G'$
in general position. Since for each~$i$,
$\ell_{(\cdot)}$ is the same number $\ell_i$
for all irreducible factors of $G_i/K_i$, we deduce from 
Remark~\ref{rmk:estimates}
that $\ell_i\dim G_i/K_i\leq2\dim G_i$ for all $i$. It follows that
\begin{eqnarray}\nonumber
  \ell_1\dim G/K &=& \sum_{i=1}^m\ell_i\dim G_i/K_i -\sum_{i=2}^m(\ell_i-\ell_1)\dim G_i/K_i \\ \label{ell1}
  &\leq& 2\dim G'-\sum_{i=2}^m(\ell_i-\ell_1)\dim G_i/K_i.
\end{eqnarray}

For each~$i$, the fixed point set of the nice involution
$\sigma_i:=g_i\sigma g_i^{-1}$ has a component of codimension 
at most $4+\dim G/K$.
Denote by $F_1$ a component of maximal dimension of
$M^{\sigma_1}\cap\cdots\cap M^{\sigma_{\ell_1}}$.
As in section~\ref{pf-simple}, $F_1$ is non-empty
and 
\begin{equation}\label{dimf1-estim}
\dim F_1 \geq\dim M-\ell_1(4+\dim G/K).
\end{equation}
Further, from Remark~\ref{rmk:estimates} we have
$\ell_m\leq 2\rk G_m+1\leq2\rk G'+1$.
We combine this inequality with
estimates~(\ref{ell1}) and~(\ref{dimf1-estim}),
and the assumption on $\dim M$, to write
\begin{eqnarray}\label{dimf1-new-estim}\nonumber
  \dim F_1&>&\alpha_G+\beta_G -\ell_1(4+\dim G/K) \\ \nonumber
  &\geq & \alpha_{G''} + \beta_{G''} + 2\dim G'+ 8\rk G'+ 4 -4\ell_1
  -\ell_1\dim G/K\\ 
  & \geq & \alpha_{G''} + \beta_{G''} +4(\ell_m-\ell_1)
  +\sum_{i=2}^m(\ell_i-\ell_1)\dim G_i/K_i\\ \nonumber
  &\geq&0.
\end{eqnarray}

Note that $\sigma$ does not centralize $G_1$.
The closure of the group generated by $\sigma_i$ for
$i=1,\ldots,\ell_1$ contains the transvection group of a totally geodesic
submanifold of $G_1/K_1\times\cdots\times G_m/K_m$
of maximal rank, so it is locally a product and contains $G_1$.
Therefore $F_1\subset M^{G_1}$.
Let $B_1$ be the component of $M^{G_1}$ that contains $F_1$.
Since $G_1$ is normalized by $G$ and $G$ is connected, $G$ acts on $B_1$.

We next claim that for all $i\geq\ell_1+1$ the totally geodesic submanifolds
$M^{\sigma_i}$ and $B_1$ intersect along a
  submanifold of dimension at least $\dim B_1-(4+\dim G/K-\dim G_1/K_1)$.
  Note first that $F_1\subset M^{\sigma_1}\cap B_1$, $B_1$ is $G$-invariant
  and $M^{\sigma_i}=g_ig_1^{-1}\cdot M^{\sigma_1}$, so
  $M^{\sigma_i}\cap B_1\neq\varnothing$. In order to estimate the
  codimension of the intersection, consider the normal space of
  $M^{\sigma_i}$ at a generic point~$q$. Since $G_1$ is a normal
  subgroup of $G$, $\nu_qM^{\sigma_i}$ splits as a sum $V_q\oplus W_q$
  where $V_q$ is the part contained in $T_q(G_1q)$ and $W_q$ is its orthogonal
  complement. Going back to the argument in the last 
paragraph of the proof of Lemma~\ref{nice},
note that along $T_q(G_1q)$ the codimension of
$M^{\sigma_i}$ is bounded by the dimension of the $(-1)$-eigenspace
of $\mathrm{Ad}_{\sigma_i}$, that is, $\dim G_1/K_1$, so
$\dim V_q\leq\dim G_1/K_1$ and similarly $\dim W_q\leq  4+\dim G/K
  -\dim G_1/K_1$. As a point $p\in M^{\sigma_i}\cap B_1$ is approached
  by generic points $q_n\in M^{\sigma_i}$, the numbers
  $\dim G_1q_n$, $\dim V_{q_n}$ and $\dim W_{q_n}$ stay constant,
  say $\dim V_{q_n}=r$ and $\dim W_{q_n}=s$,
  and (passing to a subsequence)
  $V_{q_n}$ converges to an $r$-dimensional subspace
  $V_p$ of $\nu_pM^{\sigma_i}$.
  Since $T_pB_1=(T_pM)^{G_1}$,
  %the slice representation
  %of $G_1$ at $p$ preserves $\nu_pB_1$,
  we obtain that $V_p$ is contained in $\nu_pB_1$.
  Now $\dim(\nu_pM^{\sigma_i}\cap\nu_pB_1)\geq r$. It follows that
  \begin{eqnarray}\nonumber
    \dim B_1-\dim(M^{\sigma_i}\cap B_1)&\leq&\dim T_pB_1-\dim(T_pB_1\cap T_p M^{\sigma_i})\\ \nonumber
    &=&\dim (T_pB_1+T_p M^{\sigma_i}) -\dim T_pM^{\sigma_i}\\ \label{b1}
    &=& \dim\nu_pM^{\sigma_i}-\dim(\nu_pB_1\cap\nu_p M^{\sigma_i})\\ \nonumber
    &\leq&(r+s)-r\\ \nonumber
    &\leq& 4+\dim G/K-\dim G_1/K_1.
  \end{eqnarray}

  Let $F_2$ be a component of maximal dimension of
  $F_1\cap M^{\sigma_{\ell+1}}\cap\cdots\cap M^{\sigma_{\ell_2}}$.
  By Frankel's theorem applied to~$B_1$ as ambient space,
$F_2\neq\varnothing$. In fact, using~(\ref{dimf1-new-estim})
and~(\ref{b1}), we obtain that 
    \begin{eqnarray*}
      \dim F_2&\geq&\dim F_1+\dim B_1\cap M^{\sigma_{\ell_1+1}}\cap\cdots\cap
      M^{\sigma_{\ell_2}}-\dim B_1\\
&>& \alpha_{G''}+\beta_{G''}+4(\ell_m-\ell_1)
      +\sum_{i=2}^m(\ell_i-\ell_1)\dim G_i/K_i\\
      &&\qquad-(\ell_2-\ell_1)(4+\dim G/K-\dim G_1/K_1) \\
&=&  \alpha_{G''}+\beta_{G''}+4(\ell_m-\ell_1)
      +\sum_{i=2}^m(\ell_i-\ell_1)\dim G_i/K_i\\
      &&\qquad-4(\ell_2-\ell_1)-(\ell_2-\ell_1)\sum_{i=2}^m\dim G_i/K_i\\
      &=& \alpha_{G''}+\beta_{G''} + 4(\ell_m-\ell_2)+\sum_{i=3}^m(\ell_i-\ell_2)\dim G_i/K_i\\
      &\geq&0.
    \end{eqnarray*}
    
The closure of the subgroup generated by $\sigma_i$ for $i=1,\ldots,\ell_2$ contains $G_1G_2$. Therefore $F_2\subset M^{G_1G_2}$. Let $B_2$ be the
    component of $M^{G_1G_2}$ that contains $F_2$. Note that $G$ acts
    on $B_2$. Proceeding by induction, we find a component $\tilde B=B_m$ of
    $M^{G'}$ that contains a component
    $F_m$ of maximal dimension of
    $M^{\sigma_1}\cap\cdots\cap M^{\sigma_{\ell_M}}\neq\varnothing$ of dimension
    \[ \dim\tilde B> \alpha_{G''}+\beta_{G''}. \]
    Note that $\tilde B$ is orientable and
    the action of $G''$ on $\tilde B$ satisfies the dimension
    hypothesis in the statement of the theorem, so if $\tilde B/G''$ has
    non-empty boundary, we can check whether~$(\mathcal S)$
is present or not and continue the process.

\subsection{$(\mathcal S)$ is present}
Then $G_{pr}$ is finite and there is a $G$-important point $p\in M$
such that $G_p^0$ is a central circle group. Set $G':=G_p^0$.
The fixed point set  $M^{G'}$ has codimension $2$ in $M$.
Let $\tilde B$ be the component
of $M^{G'}$ containing $p$. Then $\tilde B$ is orientable,
$G'':=G/G'$ acts on $\tilde B$ and
$\dim \tilde B=\dim M -2>\alpha_G+\beta_G-2=\alpha_{G''}+\beta_{G''}$,
so if $\tilde B/G''$ has non-empty boundary, we can 
check whether~$(\mathcal S)$
is present or not and continue the process.

\subsection{End of proof}

In any case, $G''$ is a connected Lie group with $\dim G''<\dim G$, so
the process must stop  after finitely many repeatitions
of the argument. We end up
with a component $B$ of the fixed point set of a normal subgroup 
$G^\infty$ of $G$
such that $B/G$ has empty boundary.
Let $N$ be the subgroup of $G$ consisting of all
elements that fix $B$ pointwise. It is clear that $N$ is a
(possibly disconnected) normal subgroup of $G$ of positive dimension
containing $G^\infty$,
$\dim B>\alpha_{G/N^0}+\beta_{G/N^0}$,
the action of $G/N$ on $B$ is effective and its orbit space has
    empty boundary. In particular, the principal isotropy group of
    $G/N$ on $B$ is trivial by~\cite[Lemma~3.1]{Wi2}. This proves
part~(a) and the first statement of part~(b).

    Since $\dim B/G>0$, the Frankel-Petrunin
    theorem for positively curved Alexandrov spaces~\cite[Theorem~3.2]{Pe}
    implies that $B/G$ meets each face of $M/G$.
    Since $B/G$ itself has no codimension one strata, it follows that
    $B/G$ is contained in each face of $M/G$.
    It follows that any isotropy group corresponding to a
    codimension one stratum of $M/G$ is contained in the principal
    isotropy of the action of $G$ on $B$, namely, $N$.
 This proves the second statement of part~(b) and 
part~(c)(i).

    Since $B/G$ is contained in the boundary of $M/G$, the isotropy (slice)
    representation of $N$ at a generic point $p\in B$ has orbit space
    with non-empty boundary, which is part~(c)(ii). 
Assume now $M$ is
simply-connected and let us show that the same holds for
    the isotropy representation of $N^0$ at $p$. There is
a principal isotropy group $G_{pr}$ contained in $N$. If $\dim G_{pr}>0$,
then $(G_{pr})^0\subset N^0$. This implies that the isotropy representation
of $N^0$ has non-trivial principal isotropy group and the desired result
follows from~\cite[Lemma~3.1]{Wi2}. It remains to discuss the case
in which $G_{pr}$ is finite. Recall $N$ contains all isotropy
groups corresponding to codimension one strata of $M/G$.
Owing to the simple-connectedness of $M$ and~\cite[Lemma~3.6]{GL},
there are no boundary components of $\mathbb Z_2$-type.
Now $N$ contains isotropy groups of dimensions $1$
or $3$ associated to codimension one strata of $X$, and then
    $N^0$ contain the corresponding identity components;
these groups give rise to
codimension one strata for the isotropy representation of $N^0$.
This proves~(c)(iii) and
completes the proof of Theorem~\ref{main}.

    \iffalse
    To see that the same holds for
    the isotropy representation of $N^0$, in view of~\cite[Lemma~3.1]{Wi2}
    we may assume that the principal isotropy group $H=G_{pr}$ is finite
    (recall $H^0\subset N^0)$ and of odd order (for otherwise $H\cap N^0$ contains a
    nice involution which will belong to the principal isotropy group of the isotropy representation
    of $N^0$). Now either $N$ contains isotropy groups of dimensions $1$
    or $3$ associated to codimension one strata of $X$ and then
    $N^0$ contain the corresponding identity components, or there is a
    codimension one stratum of $X$ of type $\Z_2$ whose nice involution
    belongs to $N^0$. In any case, these groups give rise to
    codimension one strata for the isotropy representation of $N^0$.
    This completes the proof of Theorem~\ref{main}.
\fi

\section{Reducible representations}\label{reducible}

The following proposition follows from~\cite[Proposition~12.1]{S},
but we provide a proof for the sake of clarity.

\begin{prop}\label{prop:reducible}
Let $\rho:G\to \OG V$ be a representation of a
compact Lie group~$G$ with orbit space $X=V/G$.
Assume $V=V_1\oplus V_2$ is a $G$-invariant decomposition,
write $\rho=\rho_1\oplus\rho_2$, denote a principal isotropy
group of $\rho_i$ by $H_i$, for $i=1$, $2$, and
put $Y_1=V_1/\rho_1(H_2)$ and $Y_2=V_2/\rho_2(H_1)$.
Then $\partial X\neq\varnothing$ if and only if
$H_2$ is non-trivial and $\partial Y_1\neq\varnothing$
or $H_1$ is non-trivial and $\partial Y_2\neq\varnothing$.
\iffalse
\ft{It appears
  from the later classification that we can replace ``or'' by ``and''
  in the last sentence,
  in case neither $V_1$ nor $V_2$ is a trivial representation, at least
  if $G$ is simple and $V_1$ and $V_2$ are irreducible.}
\fi
\end{prop}

\Pf Let $p_1\in V_1$ be a point with $G_{p_1}=H_1$.
The slice representation of $H_1$ on $\nu_{p_1}(Gp_1)$ is the sum of
a trivial component and $\rho_2|_{H_1}$. If $\partial Y_2\neq\varnothing$,
then the orbit space of the slice representation has non-empty boundary
and hence $p_1$ projects to a point in $\partial X$.

Conversely, suppose $p=p_1+p_2\in V$
is a $G$-important point, where $p_i\in V_i$.
Then the slice representation $(G_p,\nu_p:=\nu_p(Gp))$ decomposes as
the sum of a trivial component and a cohomogeneity $1$
representation. Since $\nu_p\cap V_1$ and $\nu_p\cap V_2$ are $G_p$-invariant,
this implies $G_p$ is trivial on one of them, say, $\nu_p\cap V_1=:\nu_p^1$.
We can find $p_1'\in V_1$ in the normal slice at $p_1$,
sufficiently close to $p_1$, such that $G_{p_1'}\subset
G_{p_1}$ is a principal isotropy group of $\rho_1$. By replacing $p$ by a
$G$-conjugate, we may assume $G_{p_1'}=H_1$. Put $p'=p_1'+p_2$ and
note that $G_{p'}=G_p$ and $p'$ lies in the stratum of~$p$,
since $G_p$ leaves $\nu_p^1$ pointwise fixed. In particular, $p'$ is a
$G$-important point. Moreover
\[ G_{p_1'+\lambda p_2}=(G_{p_1'})_{\lambda p_2}=(G_{p_1'})_{p_2}=G_{p'} \]
for all $\lambda\neq0$, so $p_1'+\lambda p_2$ is also a
$G$-important point, and hence $p_1'+0$ projects to $\partial X$, by
continuity. Then the slice representation of $G_{p_1'+0}$ has orbit space
with non-empty boundary, but this representation equals the
trivial action of $H_1$ on $\nu_p^1$ plus $\rho_2|_{H_1}$.
Hence $H_1\neq\{1\}$ and $\partial Y_2\neq\varnothing$. \EPf

\iffalse
\begin{cor}
  If $G$ is connected and $\rho_1$ has finite principal isotropy
  group $H_1$, then $\rho$ has empty boundary.
\end{cor}

\Pf In fact the elements of $H_1$ act on $V_2$ with fixed point
subspaces of codimensions larger than one because $G$ is connected......................................
\fi

\medskip

We will need the following lemma communicated to us
in much greater generality by the authors, see~\cite[Lemma 12.3]{KL}.
Recall that a map between metric spaces is called a \emph{submetry} if it maps
any given closed ball around a point onto the closed ball of the
same radius around the image point.
For a connected complete Riemannian manifold $M$
of positive curvature
and closed subgroups $G\subset H$ of isometries of $M$,
it is easy to see that the natural projection
$M/G\to M/H$ is a submetry between Alexandrov spaces of positive
curvature (on submetries, see also~\cite{BG,L0}).

\begin{lem}[Kapovitch-Lytchak]\label{lytchak}
For a compact Riemannian manifold $M$ of positive curvature
and closed subgroups of isometries $G\subset H$ of $M$,
consider the natural submetry $f:X=M/G\to Y=M/H$.
  If $\partial X\neq\varnothing$
  then $\partial Y \neq\varnothing$ (here we follow the usual convention
  that a point space has non-empty boundary).
\end{lem}

\Pf Suppose, to the contrary, that $\partial Y=\varnothing$.
Since $Y$ has no strata of codimension one,
by~\cite[Lemma 4.1]{LT} we can find an infinite $H$-horizontal
geodesic $\gamma$
in $M$ which meets no singular $H$-orbits and thus projects to a
geodesic $\gamma''$ in the Alexandrov space $Y$.

Let $\gamma'$ be the projection of $\gamma$ to $X$.
We can assume $\gamma$ was chosen so that $\gamma'$ starts
at a point in $X\setminus\partial X$.
Note that  $\gamma'$ is a horizontal lift of $\gamma''$ under $f$ and hence
a geodesic in the compact Alexandrov space $X$.
By positive curvature of $X$,
the distance function to $\partial X$ is
strictly concave and thus $\gamma'$ must meet
$\partial X$.

On the other hand, using~\cite[Lemma 4.1]{LT} again, we may assume
$\gamma$ was chosen so
that $\gamma'$ meets $\partial X$ at a point $x$ belonging to
a codimension one stratum.
Then $\gamma'$ is a concatenation of geodesics that satifies
the reflection law at $x$, and hence cannot be locally minimizing at~$x$,
which is a contradiction. \EPf

\begin{cor}\label{cor:reducible}
 Let $\rho:G\to \OG V$ be a representation of a
compact connected simple 
 Lie group~$G$ with no trivial components and orbit space $X=V/G$.
Assume $V=V_1\oplus V_2$ is a $G$-invariant decomposition,
write $\rho=\rho_1\oplus\rho_2$, and put $X_1=V_1/\rho_1(G)$ and $X_2=V_2/\rho_2(G)$. 
If $\partial X\neq\varnothing$ then
$\partial X_1\neq\varnothing$ and~$\partial X_2\neq\varnothing$.
\end{cor}

\Pf Let $H_1$, $H_2$, $Y_1$ and $Y_2$ be as in 
Proposition~\ref{prop:reducible};
by this proposition, say $H_1\neq\{1\}$ and $\partial Y_2\neq\varnothing$.
We claim that 
$H_1\supsetneq\ker\rho_1$. In fact, otherwise
$H_1=\ker\rho_1\supsetneq\{1\}$; this, together with the assumption 
that $G$ is simple, yields that $H_1$ is a finite subgroup,
but $G$ is connected and no element of $H_1$ can act
on~$V_2$ as a reflection on a hyperplane,
which is a contradiction to~$\partial Y_2\neq\varnothing$.
Now it follows from the claim that 
$\partial X_1\neq\varnothing$~\cite[Lemma~3.1]{Wi2}.
Finally, $\partial Y_2\neq\varnothing$ is equivalent to
$S(Y_2)=S(V_2)/\rho_2(H_1)$ having non-empty boundary.
The natural projection $S(Y_2)\to S(X_2)$ is a submetry,
so Lemma~\ref{lytchak} implies that
$\partial S(X_2)\neq\varnothing$ and
hence~$\partial X_2\neq\varnothing$. \EPf

\section{Applications}\label{appl}

\subsection{Representations of simple groups}\label{repr-simple}

\textit{Proof of Theorem~\ref{classif}.}
We only sketch the main ideas.
 The main tools
are Theorem~\ref{simple}, Lemma~\ref{nice}, Proposition~\ref{prop:reducible}
and Corollary~\ref{cor:reducible}.
We list the related invariant $\mathcal L_G$ 
for the simple Lie groups~$G$ in Table~4.
The full calculations are
too long to reproduce here and we refer the reader instead to
the unpublished manuscript~\cite{GKW2}.

\begin{table}[t]
{\small
\[ \begin{array}{|c|c|}
  \hline
  G &  \mathcal L_G \\
  \hline
  \SU 2 & 18 \\
  \SU n\ \mbox{($n\geq3$)} & 2n^2+2n \\
  \SO n\ \mbox{($n\geq3$)} & n^2+3n \\
  \SP3 & 48 \\
  \SP4 & 72 \\
  \SP5 & 102 \\
  \SP n\ \mbox{($n\geq6$)}& 4n^2 \\
    \G & 36 \\
    \F & 96 \\
    \E6 & 132 \\
    \E7 & 222 \\
    \E8 & 396 \\
    \hline
\end{array} \]
\smallskip
\begin{center}
  \sc Table~4: The invariant $\mathcal L_G$ for a compact
  connected simple Lie group~$G$.
  \end{center}
}
\end{table}

In order to obtain Table~1 (the irreducible case),
for each simple group we bound the dimension of the candidate
representations using Theorem~\ref{simple}.
To exclude the representations whose
dimension fall within that bound but are not listed
in Table~1, we check that they fail to satisfy another
necessary condition for having non-empty boundary,
namely, the existence of a nice involution (Lemma~\ref{nice}).
We refer to the Borel-de Siebenthal classification of
involutions of compact connected simple Lie
groups~\cite[Theorem~8.10.8]{Wo}. One can explicitly list
the involutions and compute
the codimensions of their fixed point sets to see
many involutions
that disobey the bound $4+\dim G/K$ in the definition of nice
involution. We also use some techniques from~\cite{GL}. 

%There is one case that resists the arguments above, which we
%deal in Lemma~\ref{g2} below.

To see that the orbit space of $(\Spin{11},\Q^{16})$ has non-empty boundary,
note that the slice representation at a highest weight vector
$(\SU5,\C^5\oplus\Lambda^2\C^5)$
(up to a trivial component of dimension~$3$)
has non-empty boundary in the orbit space.

In order to obtain Table~2 (the reducible case),
Corollary~\ref{cor:reducible} says that
we need only to check which sums of representations in Table~1
have orbit space with non-empty boundary. Here we can first apply
the dimension estimate given by Theorem~\ref{simple}, and then proceed
with the criterion given by Proposition~\ref{prop:reducible}. \EPf

\subsection{Quaternionic representations}\label{sec:quat}

\textit{Proof of Corollary~\ref{quat}.}
It is equivalent to show that the tangent spaces of the $\hat\rho(G)$-
and $\hat\rho(\SP1)$-orbits at a regular point of $\rho$ 
meet in zero only. Therefore we may assume that $G$ is a maximal closed connected subgroup of $\SP V$. According to Dynkin, $G$ is one of the following
($n=\dim_{\mathbb H}V$):
\begin{enumerate}
\item[(i)] $\U n$;
\item [(ii)] $\SP k\times\SP{n-k}\quad\mbox{($1\leq k<n$)}$;
\item[(iii)] $\SO k\otimes\SP{n/k}\quad\mbox{($3\leq k\leq n$)}$;
\item[(iv)] a simple group.
\end{enumerate}
We note that in cases~(i) and~(ii) the representation is reducible,
contrary to our assumption.

In case~(iii), the 
connected principal isotropy group of~$\hat\rho$ 
is contained in $G$, which is sufficient.
Indeed, the connected principal isotropy group of 
$\SO k\otimes\SP\ell\SP1$ is given by (cf.~\cite[p.~72]{GPo}):
\[ \left\{ \begin{array}{ll} 
           \SO{k-4\ell} & \mbox{if $k\geq4\ell+2$;} \\
           \SP{\ell-k} & \mbox{if $\ell\geq k+1$;} \\
           \{1\} & \mbox{otherwise.} 
         \end{array} \right. \]

  In case~(iv) we use Theorem~\ref{classif}. If the cohomogeneities
of $\rho$ and $\hat\rho$ do not differ by~$3$, 
  the principal isotropy group $\hat G_{pr}$  
  of $\hat\rho$ is positive-dimensional. If, in addition,
 the principal isotropy group
$G_{pr}$ of~$\rho$ is non-trivial,  
then the orbit space of $\rho$ has non-empty boundary~\cite[Lemma~3.1]{Wi2}
and $\rho$ must be listed in Table~1.
Now~$\rho$ is one of:
\begin{equation}\label{q-type-neb}
  (\Spin{11},\Q^{16}),\ (\Spin{12},\Q^{16}),\ (\SU6,\Lambda^3\C^6),\
  (\SP3,\Lambda^3_0\C^6),\ (\E7,\Q^{28}).
\end{equation}
In the first representation we have a non-maximal group, as
the half-spin representation of $\Spin{12}$ restricts to
the spin representation of $\Spin{11}$. For the remaining
four representations it is true that
$c(\rho)=c(\hat\rho)+3$ (see e.g.~\cite[Table~A]{HH}). 

There remains the case in which $\dim\hat G_{pr}>0$ and $G_{pr}$ is trivial.  
By the argument in Lemma~\ref{nice} and Remark~\ref{rmk:nice}, 
there is a nice involution $\sigma\in\hat G_{pr}$ such that
$\sigma^2=1$ and 
\begin{equation}\label{impr-nice}
\dim V-\dim V^\sigma\leq\dim\hat G/\hat G^\sigma.
\end{equation}  
Now $\sigma=(\sigma_1,\sigma_2)\in G\times\SP1$, where $\sigma_2=\pm1$.
Owing to the fact that $G_{pr}$ is trivial, $\sigma_2=-1$; 
further, $\sigma_1\neq1$ as $\sigma$ is not central in $\hat G$. 
Now $G$ acts almost effectively on the symmetric 
space of inner type $G/G^{\sigma_1}=\hat G/\hat G^\sigma$, and we 
can apply Proposition~\ref{dense} as in section~\ref{pf-simple} 
to obtain that $\dim S(V)-\dim S(V)^G\leq\ell_{G/G^{\sigma_1}}\dim G/G^{\sigma_1}
\leq \max_K\{\ell_{G/K}\dim G/K\}:=\hat{\mathcal L}_G$, where 
$K$ runs through all symmetric subgroups of $G$ with maximal rank,
and $S(V)$ denotes the unit sphere of $V$.
Due to the irreducibility of $\rho$, we have $S(V)^G=\varnothing$,
so we deduce from this inequality that $\dim V\leq\hat{\mathcal L}_G$.

The compact simple Lie groups admitting irreducible representations
of quaternionic type are listed in~\cite[p.~71]{GPo}, where also the
minimal dimension of such a representation (of cohomogeneity at least $2$)
is given. In Table~5 we list the values of $\hat{\mathcal L}_G$ for those 
groups. Running through irreducible representations of quaternionic type
of $G$ of cohomogeneity at least~$2$ and dimension at most $\hat{\mathcal L}_G$,
we precisely obtain those listed in~(\ref{q-type-neb}) and 
$(\Spin{13},\Q^{32})$. 

We finally show that~(\ref{impr-nice}) cannot hold 
for the latter representation. Indeed in this case
$V^\sigma=V^{-\sigma_1}$, and the calculation in~\cite[\S2.2]{GKW2}
shows that $\dim V^{-\sigma_1}=\frac12\dim V= 64$, so $V^\sigma$ 
has codimension $64$, which is bigger than 
$\dim\hat G/\hat G^\sigma=\dim G/G^{\sigma_1}\leq \dim\Spin{13}/(\Spin6\times\Spin7)=42$. 

\begin{table}[t]
{\small
\[ \begin{array}{|c|c|}
  \hline
  G &  \hat{\mathcal L}_G \\
  \hline
  \SU 2 & 6 \\
  \SU n\ \mbox{($n\geq3$)} & 2n^2-2n \\
  \SO n\ \mbox{($n\geq3$)} & n^2-n \\
  \SP n\ \mbox{($3\leq n\leq6$)} & 3n^2+3n \\
  \SP n\ \mbox{($n\geq7$)} & 4n^2-4n \\
  \E7 & 210 \\
    \hline
\end{array} \]
\smallskip
\begin{center}
  \sc Table~5: The invariant $\hat{\mathcal L}_G$ for some compact
  simple Lie groups~$G$.
  \end{center}
}
\end{table}
\EPf

\medskip

We will use the following lemma
in the proof of Corollary~\ref{cor:red}.

\begin{lem}\label{lem:quat}
  Let $\rho:G\to \OG V$ be an irreducible
  representation of a compact Lie group of quaternionic type
  and cohomogeneity at least two.
  Assume $\tau:H\to \OG W$ is a %minimal
  reduction of $\rho$.
  Then $\tau$ is also of quaternionic type.
\end{lem}

\Pf By assumption, the centralizer of $\rho(G)$ in $\OG V$
contains an $\SP1$-subgroup. Due to Corollary~\ref{quat},
this subgroup induces an $\SP1$- or $\SO3$-group of isometries
of $X:=V/G=W/H$.
%Since $H_{pr}$ is trivial, the natural map
%$N_{\sf O(W)}(\tau(H))\to \mathrm{Isom}(X)$ is injective.
By~\cite[Theorem~A]{mendes}, any isometry in the identity
component of the isometry group of $X$ is induced by an
element in the centralizer of $\tau(H)$ in $\OG W$. We deduce that
this centralizer has dimension at least~$3$.
Since $\tau$ is irreducible~\cite[Lemma~5.1]{GL},
this implies it is of quaternionic type. \EPf

\subsection{Representations of simple groups, continued}

\textit{Proof of Corollary~\ref{cor:red}.} A representation can admit a
non-trivial reduction only if it has non-empty boundary
in the orbit space.
Therefore, in view of Table~1, it suffices to prove that the spin
representation~$\rho$ of $G=\Spin{11}$ on $V=\Q^{16}$
admits no non-trivial reductions. For later use, recall that
its principal isotropy group is trivial and its cohomogeneity is~$9$.

Suppose, to the contrary, that $\rho$ admits non-trivial
reductions and choose a \emph{minimal} reduction
$\tau:H\to\OG W$, that is, $\tau$ satisfies $W/H=V/G=X$,
$\dim H<\dim G=55$ and $\dim H$ is as small as possible.
Then $H_{pr}$ is trivial.
%The normalizer of $G$ in $\OG V$ contains a $\SP1$-group
%which induces an effective $\SO3$-isometric action on $X$. We can lift isometries
%of $X$ to $W$ by Lemma~\ref{mendes},
Since $\rho$ is of quaternionic type,
by Lemma~\ref{lem:quat} also $\tau$ is of quaternionic type.
In particular,
$H$ is semisimple. Since $\rho$
is not toric~\cite{GL2}, it also follows from~\cite[Theorem~1.7]{GL} that
$\tau^0=\tau|_{H^0}$ is irreducible.

Next, we need to analyse irreducible representations of quaternionic type
(of dimension $<64$) of compact connected semisimple Lie groups (of dimension $<55$)
of cohomogeneity~$9$.

Assume first $H^0$ is simple. It is easy to list irreducible
representations of quaternionic type of simple groups of low dimension
and estimate their cohomogeneities. This yields $H^0=\SP1$ and
$W=\Q^3$. In this case, $W/H^0$ has empty boundary, so $H/H^0$ is generated by elements
that act on $W/H^0$ as reflections. It follows that there is an element
$\sigma\in H\setminus H^0$ of order $2$ fixing a $H$-important, $H^0$-regular
point~\cite[\S2.2, \S4.3]{GL} and
\[ \dim W-1=\dim H-\dim Z_H(\sigma) + \dim W^\sigma, \]
where $Z_H(\sigma)$ denotes the centralizer of $\sigma$ in $H$, that
is,
\[ \dim W^\sigma=8+\dim Z_H(\sigma). \]
Note that $\dim Z_H(\sigma)=1$ or $3$ is odd. Due to~\cite[Lemma~11.1]{GL},
$\dim W^\sigma$ is even, and we reach a contradiction.

We now assume $H^0$ is not simple.
We can write $H^0=H_1\times H_2$, $W=W_1\otimes_{\mathbb R}W_2$,
$\tau=\tau_1\otimes\tau_2$, where $\tau_1$ is of real type and $\tau_2$ is of
quaternionic type. It follows from~\cite[Lemma~12.1]{GL}
that the cohomogeneity
\[ c(\SO m\otimes\SP n)\geq c(\SO 3\otimes\SP 2)\geq 3\cdot 8 - (10+3)=11 \]
for $m\geq3$ and $n\geq2$ (see also~\cite[Lemma~3.5]{Go}),
so we must have $\tau_2=(\SP1,\Q)$. It follows
that $\dim W_1<16$.

Let $p_i\in W_i$ be $H_i$-regular, for $i=1$, $2$. We estimate the cohomogeneity
$c(\tau)$ by going to the slice at $p=p_1\otimes p_2$, as follows. The normal
space $\nu_p(H^0p)$ decomposes as $\nu_{p_1}(H_1p_1)\otimes\R p_2\oplus
(\nu_{p_1}(H_1p_1)\ominus\R p_1)\otimes\R^3\oplus T_{p_1}(H_1p_1)\otimes\R^3$, and
the connected $H^0$-isotropy at $p$ has the form $(H_1)_{p_1}\times\{1\}$, acting thus
  trivially on the $\R^3$-factors and on the $\nu_{p_1}(H_1p_1)$-factors.
  Therefore the cohomogeneity
\begin{eqnarray*}
  c(\tau) &=& c(\tau_1) + 3(c(\tau_1)-1) +c((H_1)_{p_1},3(T_{p_1}(H_1p_1))) \\
  &\geq&4c(\tau_1)-3+c(\SO{m_1},3\R^{m_1}) \quad \text{($m_1=\dim H_1p_1$)}\\
  &=&4c(\tau_1)+3.
\end{eqnarray*}
Now $c(\tau)=9$ implies $c(\tau_1)=1$. From the classification of transitive
linear actions on spheres, we deduce that $\tau_1$ is one of
\[ (\SO n,\R^n),\ (\G,\R^7),\ (\Spin7,\R^8),\ (\SP n\SP1,\R^{4n}); \]
the cohomogeneity of $\tau$ becomes, respectively, $\leq7$, $\geq11$, $8$
 and $\geq16$,
a contradiction. This shows that a non-trivial reduction of~$\rho$ cannot exist. \EPf

\subsection{Isometric actions of certain simple groups}\label{non-lin}

\begin{lem}\label{asys}
  Let $M=G/K$ be a connected irreducible symmetric space, where
  $G$ is the transvection group of~$M$ and $K$ is connected.
  Assume $M$ is not of Hermitian type.
  Consider the isotropy
  representation of $K$ on the tangent space at the basepoint
  and denote by $K_{pr}$ its principal isotropy group.
  Then $N_K(K_{pr})/K_{pr}$ is finite.
\end{lem}

\Pf Write $\Lg=\Lk+\Lp$ for the decomposition of the Lie algebra of
$G$ into the $\pm1$-eigenspaces of the involution. There is a
Cartan subspace $\La$ of $\Lp$ such that $K_{pr}=Z_K(\La)$.
Let $k\in N_K(K_{pr})$. The action of $k$ on $\Lp$ must
  preserve the $K_{pr}$-isotypical decomposition of $\Lp$. In
  particular, $k$ stabilizes the $K_{pr}$-fixed point set in $\Lp$.
  Since $M$ is not of Hermitian type, the latter is $\La$~\cite[p.~11]{St}.
  We get an inclusion $N_K(K_{pr})\to N_K(\La)$ inducing an injective
  homomorphism (in fact, an isomorphism)
  $N_K(K_{pr})/K_{pr}\to N_K(\La)/Z_K(\La)$, where
  the target group is finite (it is the ``little Weyl group''
  of $M$); this implies the desired result. \EPf

  \medskip

\textit{Proof of Corollary~\ref{some-simple}.}
Suppose we are given a polar action of $G$ on $M$.
Then there are singular orbits~\cite[Lemma~2.1]{FGT}.
In particular we can find $p\in M$ and a positive dimensional
isotropy group $G_p$. The slice representation at $p$ is polar
and has orbit space with non-empty boundary. It follows that $p$
projects to the boundary of $X$.

Conversely, assume $\partial X\neq\varnothing$.
Due to Theorem~\ref{simple}, $M^G$ is non-empty and $\dim M^G\geq1$;
as in the proof of Theorem~\ref{main},
it follows that any component of $M^G$ of positive dimension is contained
in $\partial X$. In particular, $G$ has a fixed point $p\in M$
and the isotropy representation $(G,T_pM)$ has orbit space
with non-empty boundary. In case $G=\SU2$, Tables~1 and~2 say that
$T_pM=\C^2$, up to a trivial representation.
Now $G$ acts transitively on the normal sphere
to the component of $M^G$ through $p$, so $M$ is fixed point
homogeneous and the result follows from~\cite[Classification
  Theorem~2.8]{GS}.

In the other cases, Tables~1 and~2  say that
the isotropy representation of $G$ on $T_pM$ is the isotropy
representation of an irreducible symmetric space, 
not of Hermitian type,
up to a trivial representation.
It follows from Lemma~\ref{asys} that the normalizer
of the principal isotropy group $G_{pr}$
is a finite extension thereof, which means the $G$-action on~$M$
is asystatic and, in particular, polar~\cite{GK}. Now we can finish
by using~\cite[Theorem~A]{FGT}. \EPf

\providecommand{\bysame}{\leavevmode\hbox to3em{\hrulefill}\thinspace}
\providecommand{\MR}{\relax\ifhmode\unskip\space\fi MR }
% \MRhref is called by the amsart/book/proc definition of \MR.
\providecommand{\MRhref}[2]{%
  \href{http://www.ams.org/mathscinet-getitem?mr=#1}{#2}
}
\providecommand{\href}[2]{#2}

%\bibliographystyle{amsalpha}
%\bibliography{ref}

\end{document}